\documentclass[english]{article} 
\usepackage[absolute,overlay]{textpos} 
\usepackage{amsthm}
\usepackage{amsmath}
\usepackage{amssymb}
\usepackage[utf8]{inputenc}
\usepackage[english]{babel}
\usepackage{graphicx}
\usepackage{geometry}
\usepackage{caption}
\usepackage{subcaption}
\usepackage{pdfpages}
\usepackage{comment}
\usepackage{bm}
\usepackage{etoolbox}
\usepackage{mathtools}
\usepackage{tikz}

\usepackage{pgfplots}
\usepackage{pgfplotstable}
\usepackage{stmaryrd} 

\usepackage{hyperref}

\usepackage[noabbrev,capitalise]{cleveref}
\crefformat{appendix}{#2#1#3}

\usetikzlibrary{fillbetween}
\usepgfplotslibrary{fillbetween}

\usetikzlibrary{shapes.misc}
\usetikzlibrary{math} 

\usetikzlibrary{arrows}

\definecolor{darkspringgreen}{rgb}{0., 0.55, 0.3}
\definecolor{dartmouthgreen}{rgb}{0.05, 0.5, 0.06}
\definecolor{etonblue}{rgb}{0.59, 0.78, 0.64}
\definecolor{airforceblue}{rgb}{0., 0.4, 0.66}
\definecolor{arylideyellow}{rgb}{0.91, 0.84, 0.42}
\definecolor{emerald}{rgb}{0.31, 0.78, 0.47}
\definecolor{uclagold}{rgb}{1.0, 0.7, 0.0}
\definecolor{cadmiumorange}{rgb}{0.93, 0.53, 0.18}

\numberwithin{theorem}{section}

\numberwithin{definition}{section}
\newtheorem{remark}{Remark}
\numberwithin{remark}{section}

\newcommand{\uvec}[2][3]{\boldsymbol{#2\mkern-#1mu}\mkern#1mu}

\newcommand\norm[1]{\left\lVert#1\right\rVert}

\newcommand{\dt}{\Delta t}

\newcommand{\by}{\uvec{y}}

\newcommand{\undu}[0]{\underline{\uvec{u}}}

\graphicspath{{./figures/}}

\begin{document}
\author{L. Micalizzi\footnote{Affiliation: Department of Mathematics, North Carolina State University, Raleigh, NC 27695, USA. Email: lmicali@ncsu.edu. (Corresponding author)} and D. Torlo\footnote{Affiliation: Dipartimento di Matematica ``Guido Castelnuovo'', Università di Roma La Sapienza, I-00185 Roma, Italy. Email: davide.torlo@uniroma1.it.}}
\title{New Efficient Implicit-Explicit Deferred Correction methods}

\maketitle

\abstract{

	In this work, we investigate implicit-explicit (IMEX) arbitrary high-order Deferred Correction (DeC) methods for the approximation of ordinary differential equations (ODEs).
	Such schemes are characterized by an iterative procedure that increases the order of accuracy by one at each iteration.
	More precisely, we study an efficient modification based on the introduction of interpolation processes between consecutive iterations, with the aim of systematically matching the accuracy achieved at each iteration with the order of the discretization employed.
	On the one hand, this modification leads to computational advantages, since the low-order iterations are performed on cheaper lower-order discretization structures; on the other hand, it endows the methods with a natural $p$-adaptive character, which is particularly appealing in the context of practical applications.
	We investigate this modification for two families of DeC schemes, providing numerical validation, efficiency assessments, and stability region plots.
	The numerical validation includes several examples involving stiff ODEs and partial differential equations (PDEs) with high-order spatial derivatives.
	The ability of the modified schemes to provide high-fidelity results at reduced computational cost, as well as the effectiveness of the adaptive strategy, is demonstrated through the numerical experiments.

}


\section{Introduction}

Many natural and technological processes can be modeled through ordinary differential equations (ODEs) and partial differential equations (PDEs).
Given the impossibility of determining exact solutions in concrete applications, several strategies have been proposed to numerically approximate them.
In such a context, high-order schemes have been proven to be particularly efficient, as they allow one to achieve smaller errors with lower computational resources.
This is the main reason for the increasing attention shown in recent years towards the Deferred Correction (DeC) framework, as it allows one to push the order of accuracy to arbitrarily high values through a systematic construction.
In fact, after its original introduction in 1949~\cite{fox1949some}, the DeC approach became popular in 2000 through the work of Dutt, Greengard and Rokhlin~\cite{Decoriginal}, where it has been effectively employed for the construction of arbitrary high order methods for ODEs.
Since then, many follow-ups and developments have been proposed, see \cite{liu2008strong,ong2020deferred,ketcheson2014comparison,christlieb2009comments,christlieb2010integral,Decremi,micalizzi2023efficient,lore_phd_thesis} for explicit methods, and \cite{minion2003semi,minion2004semi,layton2005implications,boscarino2016error,boscarino2018implicit,minion2015interweaving,speck2015multi,hamon2019multi,abgrall2020high} for implicit or implicit--explicit (IMEX) methods, as well as references therein for a non--exhaustive literature.

A distinctive feature of DeC methods is the presence of an iterative procedure that gains one order of accuracy at each iteration.
More in detail, the DeC construction is based on fixing a high order discretization of the problem under investigation and setting up an iteration process converging to its solution.
Here, in an IMEX setting, we investigate efficient modifications based on introducing interpolation processes between the iterations to systematically match the order of accuracy of the iteration structures with the one required by the specific iteration.
This brings two main advantages: the new schemes are characterized by higher efficiency, as the cost of low-order iterations is effectively reduced thanks to the employment of lower-order structures; moreover, they possess a natural adaptive character, for which the order of accuracy does not need to be fixed a priori and the iterative procedure can be performed up to a user-defined tolerance or up to the violation of user-defined criteria.
We provide a systematic assessment of their performance, stability plots, applications to time-dependent PDEs via the method-of-lines and to adaptivity.

This work constitutes the extension of~\cite{micalizzi2023new} to an IMEX setting. Similar investigations have been performed for ADER schemes in~\cite{micalizzi2023efficient,veiga2024improving}, which can in fact be interpreted as DeC methods, see~\cite{han2021dec,offner2025analysis,lore_phd_thesis}.
%
Related efficiency mechanisms have been investigated in the Spectral Deferred Correction (SDC) literature, in particular through ladder strategies for semi-implicit and multi-implicit SDC methods, where reduced-order temporal or spatial discretizations are employed during the first correction sweeps to exploit their lower formal accuracy~\cite{minion2003semi,layton2009efficiency}. Further related developments include multi-level SDC and Parallel Full Approximation Scheme in Space and Time (PFASST) approaches, where coarse space-time representations, often coupled through full approximation scheme (FAS) corrections, are used to reduce the cost of the iterative process or to enable parallelism in time~\cite{speck2015multi,hamon2019multi,minion2015interweaving}. The present work revisits this efficiency principle within the IMEX DeC framework. More precisely, we construct interpolation-based efficient variants of two DeC families, including both solution-based and right-hand-side-based interpolation strategies, and we systematically assess their efficiency and stability properties.
Furthermore, we exploit the structure of the novel schemes to design effective and robust adaptive strategies.
The novel schemes are shown to be computationally more efficient than the original versions without substantial losses in stability properties, and their adaptive versions are effectively able to automatically select the order of accuracy to match a prescribed tolerance.

The structure of this work is the following.
In Section~\ref{sec:DeC}, we introduce the original IMEX DeC methods under investigation. 
We describe their efficient modifications, along with the adaptive versions, in Section~\ref{sec:DeC_staggered} and we investigate their stability properties in Section~\ref{sec:stability}.
In Section~\ref{sec:numerical_results}, we assess the performance of the novel methods in relation to their original versions through numerical simulations.
Finally, Section~\ref{sec:conclusions} is left for conclusions and further perspectives.

\section{IMEX DeC schemes}\label{sec:DeC}
We are interested in the numerical solution of the following Cauchy problem
\begin{equation}
	\label{eq:ODE}
	\begin{cases}
		\frac{d}{dt}\uvec{u}(t) = \uvec{S}(t,\uvec{u}(t))+\uvec{N}(t,\uvec{u}(t)),\quad t\in[0,T_f], \\
		\uvec{u}(0)=\uvec{z},
	\end{cases}
\end{equation}
where $\uvec{u}(t) \in \mathbb{R}^Q$ is the unknown solution, $\uvec{z} \in \mathbb{R}^Q$ is the initial condition, and $\uvec{S},\uvec{N}: \mathbb{R}^+_0 \times \mathbb{R}^Q \to \mathbb{R}^Q$ are two given functions defining a splitting of the right-hand side of the ODE, which is assumed to satisfy the classical assumptions guaranteeing well-posedness of the problem.
Furthermore, $\uvec{N}$ and $\uvec{S}$ are assumed to describe nonstiff and stiff dynamics, respectively, to be handled differently at the numerical level.
More specifically, $\uvec{N}$ will be treated explicitly, while $\uvec{S}$ requires an implicit treatment.

In the following, we will present two DeC methods for the previous problem.
They are one-step methods, namely, given a generic time interval $[t_n,t_{n+1}]$ with time-step $\dt:=t_{n+1}-t_{n}$, they consist of recipes to compute $\uvec{u}_{n+1}\approx \uvec{u}(t_{n+1})$ starting from a known approximation $\uvec{u}_{n}\approx \uvec{u}(t_{n})$.
For what follows, it is convenient to define here some ingredients shared by both methods. 
We introduce $M+1$ subtimenodes $t^m\in[t_n,t_{n+1}]$ such that $t_n=:t^0<t^1<\dots<t^M:=t_{n+1}$.
Their number and distribution are directly related to the order of accuracy of the resulting methods.
In particular, $M+1$ equispaced subtimenodes guarantee $(M+1)$-th order of accuracy, while $M+1$ Gauss--Lobatto subtimenodes guarantee $(2M)$-th order of accuracy~\cite{micalizzi2023new}.
Herein, we consider the latter more efficient option.

For each subtimenode $t^m$, we introduce the approximation $\uvec{u}^m\approx \uvec{u}(t^m)$, which is unknown and to be determined except for the initial subtimenode $t^0:=t_n$, for which we set $\uvec{u}^0:=\uvec{u}_n$.
Although we are interested only in the approximation at the last subtimenode $t^M:=t_{n+1}$, the methods also require approximations at the intermediate ones.

We will now present the two DeC methods.

\subsection{sDeC}\label{sec:sDeC}
This method corresponds to the original DeC scheme proposed in~\cite{Decoriginal}, for more information see~\cite{minion2003semi,layton2004conservative,huang2006accelerating,speck2015multi,micalizzi2023new}. 
Given the introduction of $M+1$ subtimenodes in $[t_n,t_{n+1}]$, the scheme relies on the following integral form of the ODE
\begin{equation}
	\label{eq:sDeC_integral_form}
	\uvec{u}(t^m)=\uvec{u}(t^{m-1})+\int_{t^{m-1}}^{t^m}\left[\uvec{S}(t,\uvec{u}(t))+\uvec{N}(t,\uvec{u}(t))\right]dt,
\end{equation}
over the generic ``small'' time interval $[t^{m-1},t^m]$, hence the name sDeC.
More specifically, the previous analytical expression is discretized in a high-order fashion as follows
\begin{align}  
	\begin{split}
		\uvec{u}^{m}=\uvec{u}^{m-1}&+\Delta t  \sum_{\ell=0}^{M} \delta^m_\ell \left[\uvec{S}(t^\ell,\uvec{u}^{\ell})+\uvec{N}(t^\ell,\uvec{u}^{\ell})\right],
	\end{split}
	\label{eq:sDeC_implicit_problem}
\end{align}
where $\delta^m_\ell:=\frac{1}{\dt}\int_{t^{m-1}}^{t^m}\psi^\ell(t)dt$ are normalized coefficients, with $\left\lbrace\psi^\ell\right\rbrace_{\ell=0,\dots,M}$ being the Lagrange interpolation polynomials associated to the subtimenodes $\lbrace t^\ell \rbrace_{\ell =0,\dots,M}$.
For $m=1,\dots,M$, Equation~\eqref{eq:sDeC_implicit_problem} constitutes an algebraic system in the unknowns $\uvec{u}^{m}$ for $m=1,\dots,M$, which must be solved to get $\uvec{u}_{n+1}\approx\uvec{u}^{M}$.

With the purpose of obtaining a scheme able to properly handle the stiffness of the problem, we approximate the solution of~\eqref{eq:sDeC_implicit_problem} through the following IMEX iterative procedure, for the iteration index $p=1,\dots,P$,
\begin{align}  
	\begin{split}
		\uvec{u}^{m,(p)}=\uvec{u}^{m-1,(p)}&+\Delta t \gamma^m \left[ \uvec{S}(t^{m},\uvec{u}^{m,(p)})- \uvec{S}(t^{m},\uvec{u}^{m,(p-1)})\right]\\
		&+\Delta t \gamma^m \left[ \uvec{N}(t^{m-1},\uvec{u}^{m-1,(p)})- \uvec{N}(t^{m-1},\uvec{u}^{m-1,(p-1)})\right]\\ 
		&+\Delta t \sum_{\ell=0}^{M} \delta^m_\ell \left[\uvec{S}(t^\ell,\uvec{u}^{\ell,(p-1)})+\uvec{N}(t^\ell,\uvec{u}^{\ell,(p-1)})\right], \text{ for } m=1,\dots,M,
	\end{split}
	\label{eq:sDeC}
\end{align}
with $\gamma^m:=\frac{t^m-t^{m-1}}{\dt}$. The vector $\uvec{u}^{m,(p)}$ represents the approximation in the $m$-th subtimenode obtained at the iteration $p$, and we set $\uvec{u}^{m,(p)}:=\uvec{u}_n$ whenever $m=0$ or $p=0$.

This DeC iterative procedure gains one order of accuracy towards the solution of the implicit system at each iteration, see~\cite{Decoriginal,minion2003semi,micalizzi2023new}, thus, for a desired order $P$, we consider $M=\left \lceil \frac{P}{2}\right \rceil $, perform $P$ iterations and finally set $\uvec{u}_{n+1}:= \uvec{u}^{M,(P)}$.


\subsection{bDeC}\label{sec:bDeC}

This method can be recovered as a particular member of the DeC family discussed in~\cite{liu2008strong}, which also includes the classical sDeC formulation. In the terminology of~\cite{ketcheson2014comparison}, it corresponds to the discrete Picard iteration case, and it is referred to here as bDeC because it is based on integral formulations over the ``big'' intervals $[t^0,t^m]$.

Notable developments include the abstract DeC framework of~\cite{Decremi}, involving an explicit continuous Galerkin formulation for hyperbolic PDEs avoiding the computational cost associated with large and sparse mass matrices; positivity-preserving schemes for ODEs and shallow water models~\cite{offner2020arbitrary,ciallella2022arbitrary,ciallella2025high}; staggered conservative schemes in primitive variables~\cite{abgrall2024staggered}; adaptive schemes~\cite{micalizzi2023new}; and asymptotic-preserving schemes~\cite{abgrall2020high,abgrall2022preliminary,chertock2026new}.

Furthermore, several investigations involving arbitrary high order frameworks~\cite{ciallella2023arbitrary,micalizzi2024novel,micalizzitoro2024,micalizzi2025algorithms,micalizzi2025force} are based on this DeC formulation.

Also in this case, the definition of the method relies on an integral formulation of the ODE, but this time on the intervals $[t^{0},t^m]$:
\begin{equation}
	\label{eq:bDeC_integral_form}
	\uvec{u}(t^m)=\uvec{u}(t_n)+\int_{t^{0}}^{t^m}\left[\uvec{S}(t,\uvec{u}(t))+\uvec{N}(t,\uvec{u}(t))\right]dt.
\end{equation}
Analogously to what was previously done, we consider a high-order implicit discretization 
\begin{align}  
	\begin{split}
		\uvec{u}^{m}=\uvec{u}_n&+\Delta t  \sum_{\ell=0}^{M} \theta^m_\ell \left[\uvec{S}(t^\ell,\uvec{u}^{\ell})+\uvec{N}(t^\ell,\uvec{u}^{\ell})\right],
	\end{split}
	\label{eq:bDeC_implicit_problem}
\end{align}
with $\theta^m_\ell:=\frac{1}{\dt}\int_{t^{0}}^{t^m}\psi^\ell(t)dt$.
The IMEX DeC iteration to solve such an algebraic system takes the form, for $p=1,\dots,P$,
\begin{align}  
	\begin{split}
		\uvec{u}^{m,(p)}=\uvec{u}_n&+\Delta t  \beta^m \left[\uvec{S}(t^m,\uvec{u}^{m,(p)})-\uvec{S}(t^m,\uvec{u}^{m,(p-1)})\right]\\
		&+\Delta t  \sum_{\ell=0}^{M} \theta^m_\ell \left[\uvec{S}(t^\ell,\uvec{u}^{\ell,(p-1)})+\uvec{N}(t^\ell,\uvec{u}^{\ell,(p-1)})\right], \text{ for } m=1,\dots,M,
	\end{split}
	\label{eq:bDeC}
\end{align}
with $\beta^m:=\frac{t^m-t^0}{\dt}$. Also here, we set $\uvec{u}^{m,(p)}:=\uvec{u}_n$ whenever $m=0$ or $p=0$. 
Concerning the accuracy and the number of iterations, the same considerations as for sDeC apply.
Note that here the difference of the explicit terms approximated with the first order approximations is not present, unlike in \eqref{eq:sDeC}, as they are approximated at $t^0$ and, since they coincide, they cancel out. 

\section{Efficient IMEX DeC schemes}\label{sec:DeC_staggered}
Here, we discuss efficient modifications of the previously presented IMEX DeC schemes. The underlying idea is to select the discretization structures, i.e., the number of subtimenodes, according to the specific accuracy achieved in the iteration.
More specifically, while in the original methods the number of subtimenodes is fixed and kept constant throughout the whole iterative process, here we start with the minimal number of subtimenodes and we alternate iterations and interpolation processes to reach the final accuracy.
Related ideas were already sketched in~\cite{minion2003semi} in the context of ladder methods, where the use of fewer subtimenodes in the lower-order iterations was considered. They were later developed more systematically in~\cite{layton2009efficiency}. In the DeC framework, interpolation-based efficient variants were investigated in the explicit setting in~\cite{micalizzi2023new}. Here, we extend this strategy to the IMEX setting and apply it to both sDeC and bDeC formulations.
As in~\cite{micalizzi2023new}, together with the standard approach in which the unknown quantity is interpolated, we also investigate an alternative approach in which the interpolation is applied to the ODE right-hand side. This choice decreases the number of required right-hand side evaluations.
We keep our notation consistent with~\cite{micalizzi2023new} and denote the modified schemes by direct reference to the interpolated quantities: ``u'' for interpolation of the solution; ``du'' for interpolation of the right-hand side (i.e., the time derivative of the solution).

The modified schemes are constructed as follows.
For order $P$, we fix the final number of subtimenodes to be $M=\left \lceil \frac{P}{2}\right \rceil $.
The subtimenodes are iteration-dependent, hence, we define the vectors $\underline{t}^{(p)}:=\left(t^{0,(p)}, \dots, t^{M^{(p)},(p)} \right)^T$ of the subtimenodes, in which we obtain the solution approximations at the $p$-th iteration. In particular, we have
\begin{align}
	M^{(p)}:=\begin{cases}
		1, \quad & \text{if}\quad p=0,\\
		p, \quad & \text{if}\quad p<M,\\
		M, \quad & \text{if}\quad M \leq p \leq P,
	\end{cases}
\end{align}
and the subtimenodes have the chosen distribution, i.e., Gauss-Lobatto. Therefore, for each $p$, $t^{0,(p)}:=t_n$ and $t^{M^{(p)},(p)}:=t_{n+1}$.

We start with $\underline{\uvec{U}}^{(0)}:=\left(\uvec{u}_n, \uvec{u}_n \right)^T$ corresponding to two subtimenodes $\underline{t}^{(0)}:=\left(t^{0,(0)}, t^{1,(0)} \right)^T$, and we perform the first iteration to get $\underline{\uvec{U}}^{(1)}=\left(\uvec{u}_n, \uvec{u}^{1,(1)} \right)^T$. 
We thus perform the second iteration making use of interpolated quantities, either solution or right-hand side, in the new subtimenodes $\underline{t}^{(2)}:=\left(t^{0,(2)},t^{1,(2)}, t^{2,(2)}\right)^T$ to get $\underline{\uvec{U}}^{(2)}=\left(\uvec{u}_n, \uvec{u}^{1,(2)},\uvec{u}^{2,(2)} \right)^T$.
We keep performing such interpolation-based iterations, for $ 2 \leq p\leq M$ and $m=1,\dots,M^{(p)}$, defined as follows
\begin{itemize}
	\item \textbf{sDeCu}\\
		\begin{align}  
			\begin{split}
				\uvec{u}^{m,(p)}=\uvec{u}^{m-1,(p)}&+\Delta t \gamma^{m,(p)} \left[ \uvec{S}(t^{m,(p)},\uvec{u}^{m,(p)})- \uvec{S}(t^{m,(p)},\uvec{u}^{*m,(p-1)})\right]\\
				&+\Delta t \gamma^{m,(p)} \left[ \uvec{N}(t^{m-1,(p)},\uvec{u}^{m-1,(p)})- \uvec{N}(t^{m-1,(p)},\uvec{u}^{*m-1,(p-1)})\right]\\ 
				&+\Delta t \sum_{\ell=0}^{M^{(p)}} \delta^{m,(p)}_\ell \left[\uvec{S}(t^{\ell,(p)},\uvec{u}^{*\ell,(p-1)})+\uvec{N}(t^{\ell,(p)},\uvec{u}^{*\ell,(p-1)})\right],
			\end{split}
			\label{eq:sDeCu}
		\end{align}
	\item \textbf{sDeCdu}\\
		\begin{align}  
			\begin{split}
				\uvec{u}^{m,(p)}=\uvec{u}^{m-1,(p)}&+\Delta t \gamma^{m,(p)} \left[ \uvec{S}(t^{m,(p)},\uvec{u}^{m,(p)})- \uvec{S}^{*m,(p-1)}\right]\\
				&+\Delta t \gamma^{m,(p)} \left[ \uvec{N}(t^{m-1,(p)},\uvec{u}^{m-1,(p)})- \uvec{N}^{*m-1,(p-1)}\right]\\ 
				&+\Delta t \sum_{\ell=0}^{M^{(p)}} \delta^{m,(p)}_\ell \left[\uvec{S}^{*\ell,(p-1)}+\uvec{N}^{*\ell,(p-1)}\right],
			\end{split}
			\label{eq:sDeCdu}
		\end{align}
	\item \textbf{bDeCu}\\
		\begin{align}  
			\begin{split}
				\uvec{u}^{m,(p)}=\uvec{u}_n&+\Delta t  \beta^{m,(p)} \left[\uvec{S}(t^{m,(p)},\uvec{u}^{m,(p)})-\uvec{S}(t^{m,(p)},\uvec{u}^{*m,(p-1)})\right]\\
				&+\Delta t  \sum_{\ell=0}^{M^{(p)}} \theta^{m,(p)}_\ell \left[\uvec{S}(t^{\ell,(p)},\uvec{u}^{*\ell,(p-1)})+\uvec{N}(t^{\ell,(p)},\uvec{u}^{*\ell,(p-1)})\right],
			\end{split}
			\label{eq:bDeCu}
		\end{align}
	\item \textbf{bDeCdu}\\
		\begin{align}  
	\begin{split}
		\uvec{u}^{m,(p)}=\uvec{u}_n&+\Delta t  \beta^{m,(p)} \left[\uvec{S}(t^{m,(p)},\uvec{u}^{m,(p)})-\uvec{S}^{*m,(p-1)}\right]\\
		&+\Delta t  \sum_{\ell=0}^{M^{(p)}} \theta^{m,(p)}_\ell \left[\uvec{S}^{*\ell,(p-1)}+\uvec{N}^{*\ell,(p-1)}\right],
	\end{split}
	\label{eq:bDeCdu}
\end{align}	
\end{itemize}
where the quantities labeled with $*$ are interpolated from the subtimenodes
$\underline{t}^{(p-1)}=\Big(t^{0,(p-1)},\allowbreak \dots,\allowbreak 
t^{M^{(p-1)},(p-1)}\Big)^T$
to the subtimenodes
$\underline{t}^{(p)}=\left(t^{0,(p)},\dots,
t^{M^{(p)},(p)}\right)^T$.
In particular, in the du variants, the quantities
$\uvec{S}^{*\ell,(p-1)}$ and $\uvec{N}^{*\ell,(p-1)}$
denote the interpolated values of
$\uvec{S}(t,\uvec{u}(t))$ and $\uvec{N}(t,\uvec{u}(t))$
at the subtimenode $t^{\ell,(p)}$, reconstructed from their values on
$\underline{t}^{(p-1)}$.
The vector $\underline{\uvec{U}}^{(M)}=\left(\uvec{u}_n, \uvec{u}^{1,(M)},\dots,\uvec{u}^{M,(M)} \right)^T$ obtained at the $M$-th iteration is associated with the $M+1$ final subtimenodes, hence, for $p=M+1,\dots,P$, we continue the iterative process on the same set of subtimenodes, without further interpolation.
Indeed, the normalized coefficients $\gamma^{m,(p)}$, $\delta^{m,(p)}_\ell$, $\beta^{m,(p)}$ and $\theta^{m,(p)}_\ell$ are associated with the subtimenodes $\underline{t}^{(p)}$.

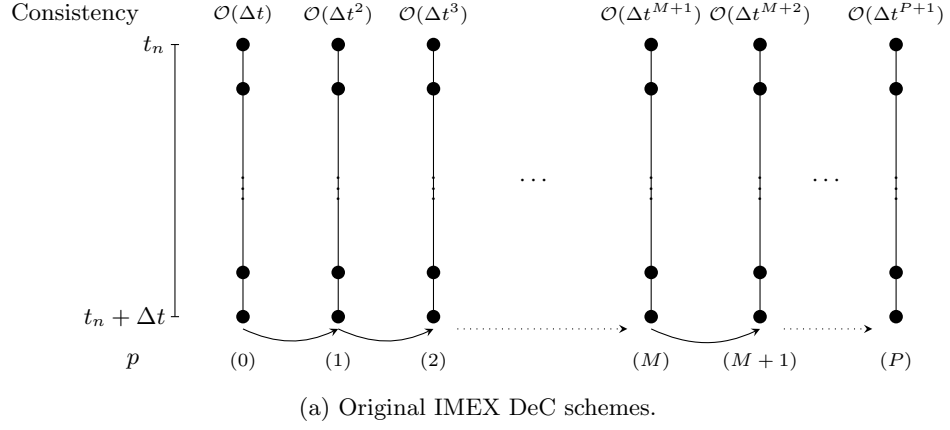
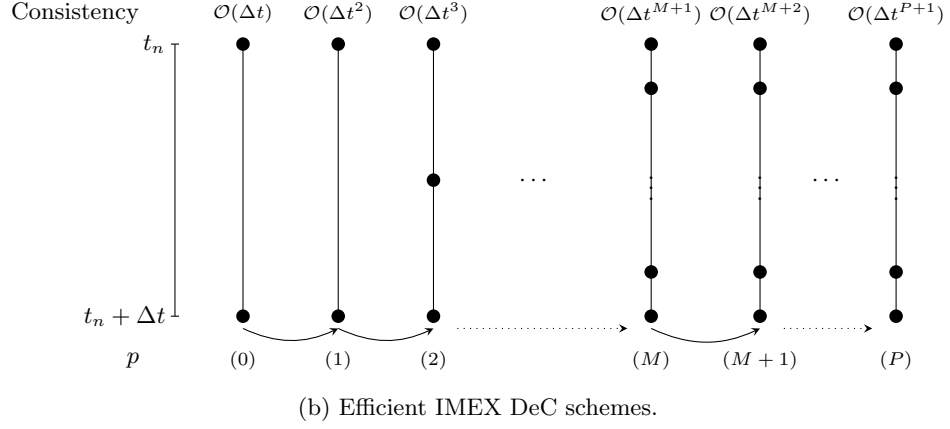
\begin{figure}
	\centering
	
	\begin{subfigure}{\textwidth}
		\centering
		\begin{tikzpicture}[
			scale=0.9,
			every node/.style={font=\small},
			dot/.style={circle, fill=black, inner sep=1.8pt},
			arrow/.style={->, >=stealth, bend right=25},
			]
			
			\def\ytop{4.0}
			\def\ybot{0.0}
			
			\def\xaxis{0.0}
			
			\def\xA{1.0}
			\def\xB{2.4}
			\def\xC{3.8}
			
			\def\xF{7.0}
			\def\xG{8.6}
			\def\xH{10.6}      
			\def\xGHmid{9.6}   
			
			\node[left] at (\xaxis,\ytop) {$t_n$};
			\node[left] at (\xaxis,\ybot) {$t_n+\Delta t$};
			
			\draw (\xaxis,\ybot) -- (\xaxis,\ytop);
			\draw (\xaxis-0.06,\ytop) -- (\xaxis+0.06,\ytop);
			\draw (\xaxis-0.06,\ybot) -- (\xaxis+0.06,\ybot);
			
			\node[left=0.35cm] at (\xaxis,\ytop+0.45) {Consistency};
			
			\node at (\xA,\ytop+0.45) {\scriptsize $\mathcal{O}(\Delta t)$};
			\node at (\xB,\ytop+0.45) {\scriptsize $\mathcal{O}(\Delta t^2)$};
			\node at (\xC,\ytop+0.45) {\scriptsize $\mathcal{O}(\Delta t^3)$};
			
			\node at (\xF,\ytop+0.45) {\scriptsize $\mathcal{O}(\Delta t^{M+1})$};
			\node at (\xG,\ytop+0.45) {\scriptsize $\mathcal{O}(\Delta t^{M+2})$};
			\node at (\xH,\ytop+0.45) {\scriptsize $\mathcal{O}(\Delta t^{P+1})$};
			
			\foreach \x in {\xA,\xB,\xC,\xF,\xG,\xH} {
				\draw (\x,\ybot) -- (\x,\ytop);
			}
			
			
			\node[dot] at (\xA,\ytop) {};
			\node[dot] at (\xA,3.35) {};
			\node at (\xA,2.0) {$\vdots$};
			\node[dot] at (\xA,0.65) {};
			\node[dot] at (\xA,\ybot) {};
			
			\node[dot] at (\xB,\ytop) {};
			\node[dot] at (\xB,3.35) {};
			\node at (\xB,2.0) {$\vdots$};
			\node[dot] at (\xB,0.65) {};
			\node[dot] at (\xB,\ybot) {};
			
			\node[dot] at (\xC,\ytop) {};
			\node[dot] at (\xC,3.35) {};
			\node at (\xC,2.0) {$\vdots$};
			\node[dot] at (\xC,0.65) {};
			\node[dot] at (\xC,\ybot) {};
			
			\node at (5.3,2.0) {$\cdots$};
			\node at (\xGHmid,2.0) {$\cdots$};
			
			
			\node[dot] at (\xF,\ytop) {};
			\node[dot] at (\xF,3.35) {};
			\node at (\xF,2.0) {$\vdots$};
			\node[dot] at (\xF,0.65) {};
			\node[dot] at (\xF,\ybot) {};
			
			\node[dot] at (\xG,\ytop) {};
			\node[dot] at (\xG,3.35) {};
			\node at (\xG,2.0) {$\vdots$};
			\node[dot] at (\xG,0.65) {};
			\node[dot] at (\xG,\ybot) {};
			
			\node[dot] at (\xH,\ytop) {};
			\node[dot] at (\xH,3.35) {};
			\node at (\xH,2.0) {$\vdots$};
			\node[dot] at (\xH,0.65) {};
			\node[dot] at (\xH,\ybot) {};
			
			\draw[arrow] (\xA,\ybot-0.18) to (\xB,\ybot-0.18);
			\draw[arrow] (\xB,\ybot-0.18) to (\xC,\ybot-0.18);
			
			\draw[dotted, ->, >=stealth] (\xC+0.35,\ybot-0.18) -- (\xF-0.35,\ybot-0.18);
			
			\draw[arrow] (\xF,\ybot-0.18) to (\xG,\ybot-0.18);
			\draw[dotted, ->, >=stealth] (\xG+0.35,\ybot-0.18) -- (\xH-0.35,\ybot-0.18);
			
			\node at (\xA,\ybot-0.65) {\scriptsize $(0)$};
			\node at (\xB,\ybot-0.65) {\scriptsize $(1)$};
			\node at (\xC,\ybot-0.65) {\scriptsize $(2)$};
			
			\node at (\xF,\ybot-0.65) {\scriptsize $(M)$};
			\node at (\xG,\ybot-0.65) {\scriptsize $(M+1)$};
			\node at (\xH,\ybot-0.65) {\scriptsize $(P)$};
			
			\node[left=0.35cm] at (\xaxis,\ybot-0.65) {$p$};
			
		\end{tikzpicture}
		\caption{Original IMEX DeC schemes.}
		\label{fig:sub_original_imex_dec}
	\end{subfigure}
	
	\vspace{1.2em}
	
	\begin{subfigure}{\textwidth}
		\centering
		\begin{tikzpicture}[
			scale=0.9,
			every node/.style={font=\small},
			dot/.style={circle, fill=black, inner sep=1.8pt},
			arrow/.style={->, >=stealth, bend right=25},
			]
			
			\def\ytop{4.0}
			\def\ybot{0.0}
			
			\def\xaxis{0.0}
			
			\def\xA{1.0}
			\def\xB{2.4}
			\def\xC{3.8}
			
			\def\xF{7.0}
			\def\xG{8.6}
			\def\xH{10.6}      
			\def\xGHmid{9.6}   
			
			\node[left] at (\xaxis,\ytop) {$t_n$};
			\node[left] at (\xaxis,\ybot) {$t_n+\Delta t$};
			
			\draw (\xaxis,\ybot) -- (\xaxis,\ytop);
			\draw (\xaxis-0.06,\ytop) -- (\xaxis+0.06,\ytop);
			\draw (\xaxis-0.06,\ybot) -- (\xaxis+0.06,\ybot);
			
			\node[left=0.35cm] at (\xaxis,\ytop+0.45) {Consistency};
			
			\node at (\xA,\ytop+0.45) {\scriptsize $\mathcal{O}(\Delta t)$};
			\node at (\xB,\ytop+0.45) {\scriptsize $\mathcal{O}(\Delta t^2)$};
			\node at (\xC,\ytop+0.45) {\scriptsize $\mathcal{O}(\Delta t^3)$};
			
			\node at (\xF,\ytop+0.45) {\scriptsize $\mathcal{O}(\Delta t^{M+1})$};
			\node at (\xG,\ytop+0.45) {\scriptsize $\mathcal{O}(\Delta t^{M+2})$};
			\node at (\xH,\ytop+0.45) {\scriptsize $\mathcal{O}(\Delta t^{P+1})$};
			
			\foreach \x in {\xA,\xB,\xC,\xF,\xG,\xH} {
				\draw (\x,\ybot) -- (\x,\ytop);
			}
			
			
			\node[dot] at (\xA,\ytop) {};
			\node[dot] at (\xA,\ybot) {};
			
			\node[dot] at (\xB,\ytop) {};
			\node[dot] at (\xB,\ybot) {};
			
			\node[dot] at (\xC,\ytop) {};
			\node[dot] at (\xC,2.0) {};
			\node[dot] at (\xC,\ybot) {};
			
			\node at (5.3,2.0) {$\cdots$};
			\node at (\xGHmid,2.0) {$\cdots$};
			
			
			\node[dot] at (\xF,\ytop) {};
			\node[dot] at (\xF,3.35) {};
			\node at (\xF,2.0) {$\vdots$};
			\node[dot] at (\xF,0.65) {};
			\node[dot] at (\xF,\ybot) {};
			
			\node[dot] at (\xG,\ytop) {};
			\node[dot] at (\xG,3.35) {};
			\node at (\xG,2.0) {$\vdots$};
			\node[dot] at (\xG,0.65) {};
			\node[dot] at (\xG,\ybot) {};
			
			\node[dot] at (\xH,\ytop) {};
			\node[dot] at (\xH,3.35) {};
			\node at (\xH,2.0) {$\vdots$};
			\node[dot] at (\xH,0.65) {};
			\node[dot] at (\xH,\ybot) {};
			
			\draw[arrow] (\xA,\ybot-0.18) to (\xB,\ybot-0.18);
			\draw[arrow] (\xB,\ybot-0.18) to (\xC,\ybot-0.18);
			
			\draw[dotted, ->, >=stealth] (\xC+0.35,\ybot-0.18) -- (\xF-0.35,\ybot-0.18);
			
			\draw[arrow] (\xF,\ybot-0.18) to (\xG,\ybot-0.18);
			\draw[dotted, ->, >=stealth] (\xG+0.35,\ybot-0.18) -- (\xH-0.35,\ybot-0.18);
			
			\node at (\xA,\ybot-0.65) {\scriptsize $(0)$};
			\node at (\xB,\ybot-0.65) {\scriptsize $(1)$};
			\node at (\xC,\ybot-0.65) {\scriptsize $(2)$};
			
			\node at (\xF,\ybot-0.65) {\scriptsize $(M)$};
			\node at (\xG,\ybot-0.65) {\scriptsize $(M+1)$};
			\node at (\xH,\ybot-0.65) {\scriptsize $(P)$};
			
			\node[left=0.35cm] at (\xaxis,\ybot-0.65) {$p$};
			
		\end{tikzpicture}
		\caption{Efficient IMEX DeC schemes.}
		\label{fig:sub_new_imex_dec}
	\end{subfigure}
	
	\caption{Progression of subtimenodes and consistency error for original and new IMEX DeC schemes.}
	\label{fig:subtimenodes_progression_imex_dec}
\end{figure}

A useful sketch of the subtimenodes progression along the iterative procedure is reported in Figure~\ref{fig:subtimenodes_progression_imex_dec}.

\begin{remark}
	According to the described algorithm, the iterations $p=1$ and $p> M$ take place without interpolation.
	Concerning $p=1$, two subtimenodes guarantee an $\mathcal{O}(\dt^2)$-interpolation error, which is consistent with the consistency error realized only after the first iteration.
	Concerning $p> M$, these iterations are meant to reach the final accuracy without adding new subtimenodes. Indeed, one could add new subtimenodes, but the resulting methods would be suboptimal.
\end{remark}

\begin{remark}
It is worth stressing that, in view of the previous remark, the modified schemes would retain higher computational advantages with respect to the original versions when equispaced subtimenodes are considered, as shown in the explicit setting~\cite{micalizzi2023new}. Here, instead, aiming at future real-world applications, we focus on Gauss-Lobatto subtimenodes, which are characterized by smaller computational costs.
\end{remark}

\begin{remark}
	The sDeCu version is closely related to the ladder strategies investigated in~\cite{minion2003semi,layton2009efficiency}, and we keep it as a reference.
\end{remark}


\subsection{Adaptive versions}\label{sec:adaptivity}
It is possible to suitably employ the new methods to design efficient $p$-adaptive schemes.
The idea, introduced in~\cite{micalizzi2023new}, consists in not fixing the final number of subtimenodes (and hence the order of accuracy) a priori.
Rather, we keep adding subtimenodes throughout the iterative process until a stopping criterion is satisfied.
Here, we consider the following criterion
\begin{equation}
	\frac{\norm{ \underline{\uvec{u}}^{M^{(p)},(p)} - \underline{\uvec{u}}^{M^{(p-1)},(p-1)}}_2}{\norm{\underline{\uvec{u}}^{M^{(p)},(p)}}_2}\leq \varepsilon,
	\label{eq:tol}
\end{equation}
where $\varepsilon$ is a user-defined convergence tolerance.
The resulting schemes are able to efficiently and automatically select the order of accuracy according to the prescribed tolerance.

\section{Linear stability}\label{sec:stability}
In this section, we investigate the linear stability of the proposed methods.
To this end, we start by observing that DeC methods for ODEs can be, in general, written as Runge--Kutta methods, see~\cite{liu2008strong,ketcheson2014comparison,han2021dec,offner2025analysis,micalizzi2023new,lore_phd_thesis}.
A generic IMEX Runge--Kutta scheme with $S$ stages on the interval $[t_n,t_{n+1}]$  reads
\begin{equation}\label{eq:RK}
	\begin{cases}
		\by^s = \uvec{u}_n + \dt \sum\limits_{r = 1}^{S} a^{\text{EX}}_{s,r} \uvec{N}(t_n+c^{\text{EX}}_r \dt, \by^r)+ \dt \sum\limits_{r = 1}^{S} a^{\text{IM}}_{s,r} \uvec{S}(t_n+c^{\text{IM}}_r \dt, \by^r), \quad s=1,\dots,S,\\
		\uvec{u}_{n+1}  = \uvec{u}_n + \dt \sum\limits_{r = 1}^{S} b^{\text{EX}}_{r} \uvec{N}(t_n+c^{\text{EX}}_r \dt, \by^r) + \dt \sum\limits_{r = 1}^{S} b^{\text{IM}}_{r} \uvec{S}(t_n+c^{\text{IM}}_r \dt, \by^r).
	\end{cases}	 
\end{equation}
The coefficients $a_{s,r}^{\text{EX}}$, $c_r^{\text{EX}}$, $b_r^{\text{EX}}$, $a_{s,r}^{\text{IM}}$, $c_r^{\text{IM}}$, and  $b_r^{\text{IM}}$ are often stored in related vector- and matrix-structures constituting the Butcher tableaux
$$\begin{array}{c|c}
	\uvec{c}^{\text{EX}}& A^{\text{EX}}  \\
	\hline
	& \uvec{b}^{\text{EX}} 
\end{array}, \quad \begin{array}{c|c}
\uvec{c}^{\text{IM}}& A^{\text{IM}}  \\
\hline
& \uvec{b}^{\text{IM}} 
\end{array},$$
which characterize the method. In particular, the matrix $A^{\text{EX}}$ is strictly lower triangular, so that the handling of $\uvec{N}$ is always explicit. Furthermore, the DeC schemes considered here are diagonally implicit Runge-Kutta (DIRK) methods, as they involve the solution of a single nonlinear system per stage (i.e., per subtimenode in each iteration), and thus the corresponding matrix $A^{\text{IM}}$ is lower triangular with nonzero diagonal entries. 
For the sake of compactness, we do not explicitly report the Butcher tableaux associated with the schemes under investigation, but only the corresponding stability results. The procedure for deriving the Butcher tableaux of the IMEX DeC methods is analogous to that described in~\cite{micalizzi2023new} for the explicit case and consists of grouping the updates performed at each iteration into block-structures that form the tableaux.

The linear stability of Runge--Kutta methods is studied on Dahlquist's equation
\begin{align}
	\begin{cases}
		\frac{d}{dt}u=\lambda u(t),\\
		u(0)=1,	
	\end{cases}\quad t\in \mathbb{R}^+_0,
\end{align}
which represents the model equation for the evolution of perturbations, with $\lambda\in\mathbb{C}$ being a complex number with negative real part, $Re(\lambda)<0.$ Since the exact solution $u=e^{\lambda t}$ is such that $\lim\limits_{t\rightarrow +\infty}u = 0$, we are interested in assessing which conditions guarantee $\lim\limits_{n\rightarrow +\infty}u_n=0$ at the discrete level.

It is worth remarking that, since the problem is linear and the interpolation is a linear operator, there is no difference in the ``u'' and ``du'' approaches in this context, as in the explicit case~\cite{micalizzi2023new}. Indeed, the interpolation and the evolution operator commute. Therefore, for stability purposes, bDeCu and bDeCdu coincide as well as sDeCu and sDeCdu.

We will consider two different approaches to the study of stability.
In Section~\ref{sec:stability_Minion}, we adopt the IMEX approach from~\cite{minion2003semi} in which a splitting of the right-hand side is considered; instead, in~\ref{sec:stability_fully_implicit}, we perform a classical stability study of the implicit part of the schemes only.
We remark here that some notable results from the explicit case do not extend to the IMEX setting.
While in the explicit case the stability regions of any bDeC, bDeCu and bDeCdu method of order $P$ coincide independently of the distribution of the subtimenodes~\cite{micalizzi2023new}, some slight differences occur in the implicit and IMEX setting.
On the other hand, the following stability results demonstrate how the efficient modifications do not negatively affect stability, as the new schemes have stability regions similar to the ones of the original methods, meaning that the gained computational efficiency and adaptive properties do not have negative drawbacks in terms of stability properties.

\subsection{Minion's stability}\label{sec:stability_Minion}
Studying the stability of an IMEX Runge--Kutta method through a classical analysis is not straightforward, as the possible complex coefficients of the implicit and explicit parts make the stability function be defined on $\mathbb C^2$, which renders the visualization of the stability region impossible.
A simplification was proposed in \cite{minion2003semi}, which consists in assuming $\lambda:=\lambda^{\text{IM}}+i\lambda^{\text{EX}}$, with $\lambda^{\text{IM}},\,\lambda^{\text{EX}} \in \mathbb R$. Then, one handles $\lambda^{\text{IM}}u$ implicitly and $i\lambda^{\text{EX}}u$ explicitly, respectively.
The rationale behind this choice is related to the fact that, usually, in PDE problems the stiff part is given by some diffusion operator involving second derivatives, whose discretizations have large real negative eigenvalues; at the same time, a simple central advection discretization has purely imaginary eigenvalues.
Thanks to the linearity of the schemes and of the problem, it is possible to express the solution in a generic iteration as 
$u_{n+1}=R_{\text{M}}(z^{\text{IM}},iz^{\text{EX}})u_{n}$, with $R_{\text{M}}(\cdot,\cdot)$ being the stability function, $z^{\text{IM}}:=\dt \lambda^{\text{IM}}\in \mathbb R$ and $z^{\text{EX}}:=\dt \lambda^{\text{EX}}\in \mathbb R$.
In particular, simple direct computations yield
\begin{equation}\label{eq:Minion_stability_function}
	R_{\text{M}}(z^{\text{IM}},iz^{\text{EX}}) = 1+ \left[z^{\text{IM}}\left(\uvec{b}^{\text{IM}}\right)^\top+iz^{\text{EX}}\left(\uvec{b}^{\text{EX}}\right)^\top\right] \left(I-z^{\text{IM}}A^{\text{IM}}-iz^{\text{EX}}A^{\text{EX}}\right)^{-1} \mathbf{1},
\end{equation}
where $I\in\mathbb R^{S\times S}$ is the identity matrix and $\uvec{1}\in \mathbb R^{S}$ is a vector of ones.
We define Minion's stability region \cite{minion2003semi} of this approach as 
\begin{equation}
	\mathcal{S}_{\text{M}}:=\left\lbrace z^{\text{IM}}+i z^{\text{EX}}\in \mathbb{C} \quad \text{s.t.} \quad \vert R_{\text{M}}(z^{\text{IM}},iz^{\text{EX}}) \vert <1 \right\rbrace.
\end{equation}

\begin{figure}
	\begin{center}
		IMEX bDeC (continuous), bDeCu/du (dashed)\\
		\includegraphics[height=0.36\textheight]{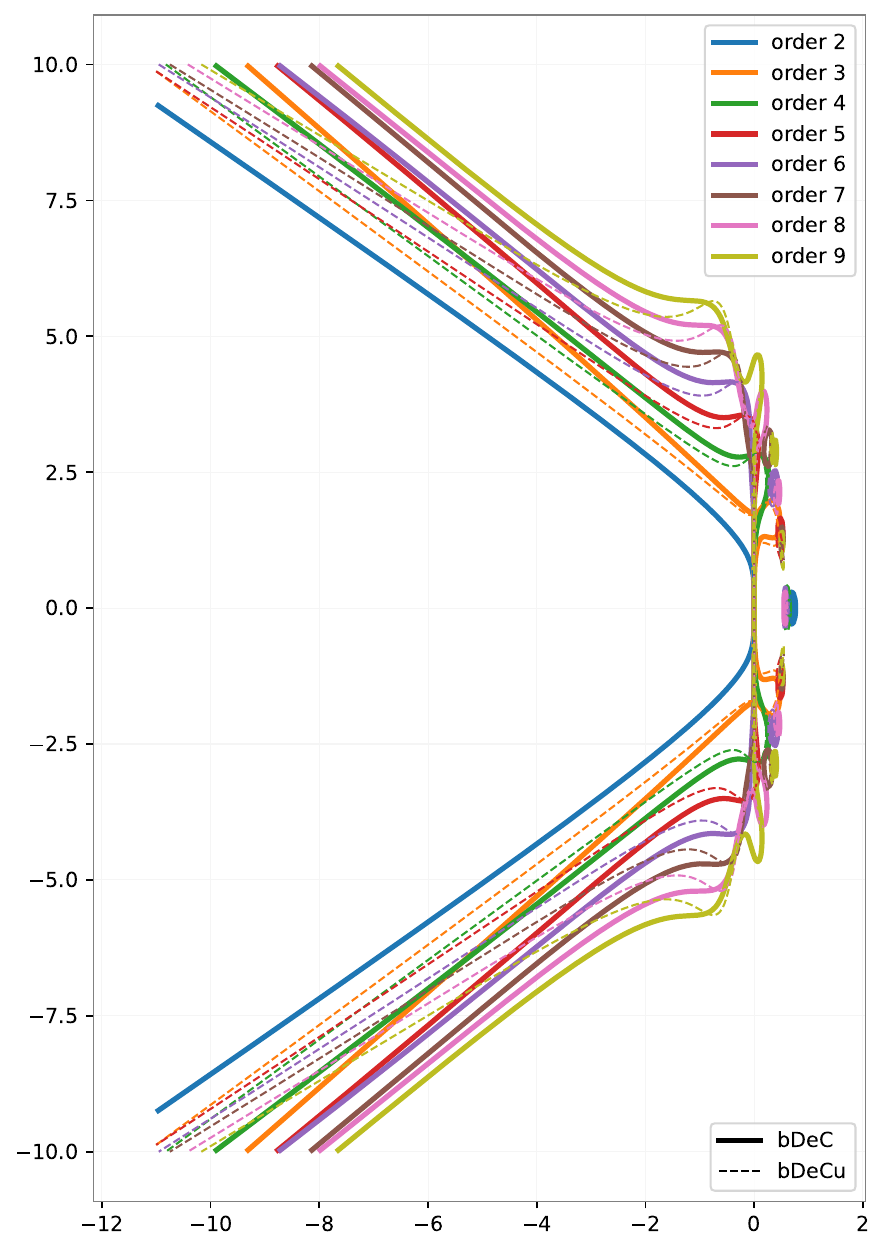}
		\includegraphics[height=0.36\textheight]{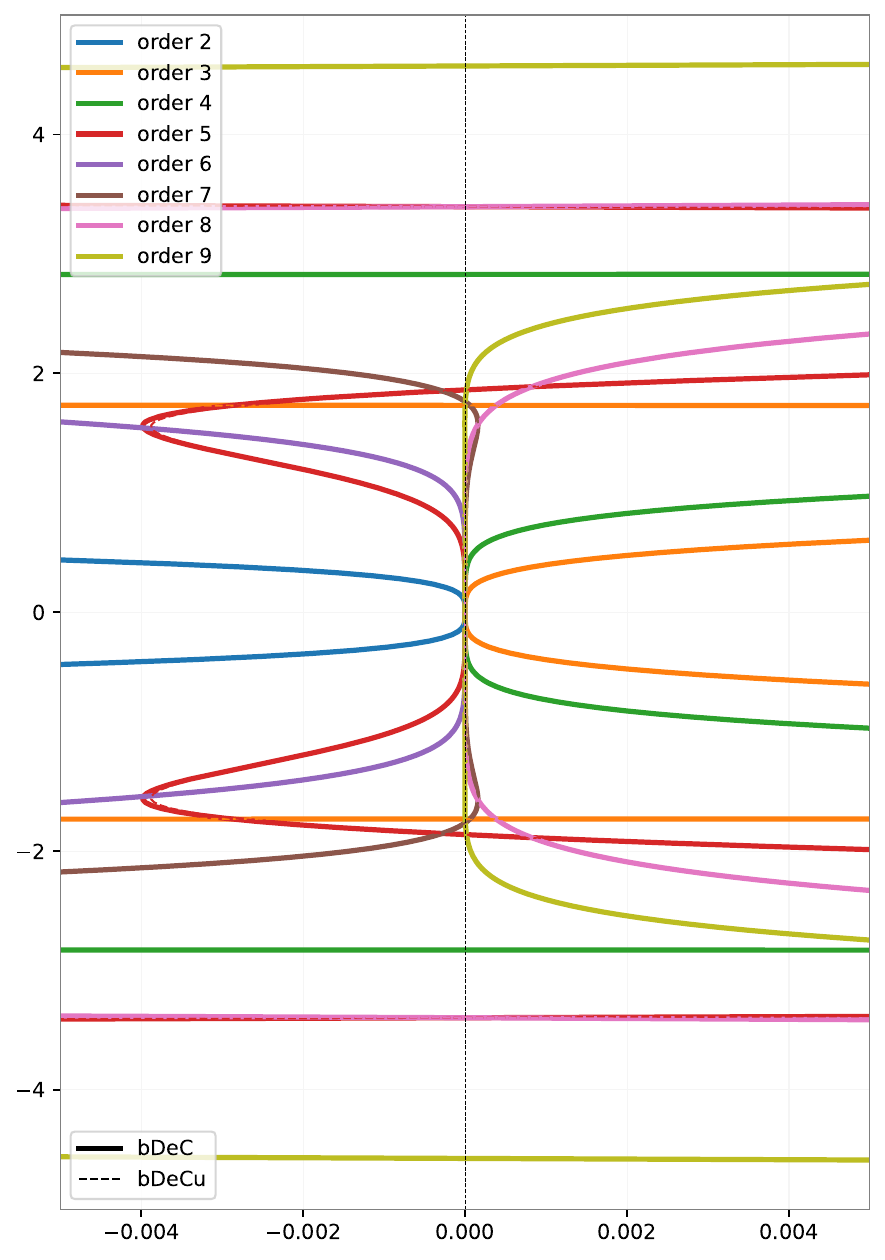}\\
		IMEX sDeC (continuous), sDeC/du (dashed)\\
		\includegraphics[height=0.36\textheight]{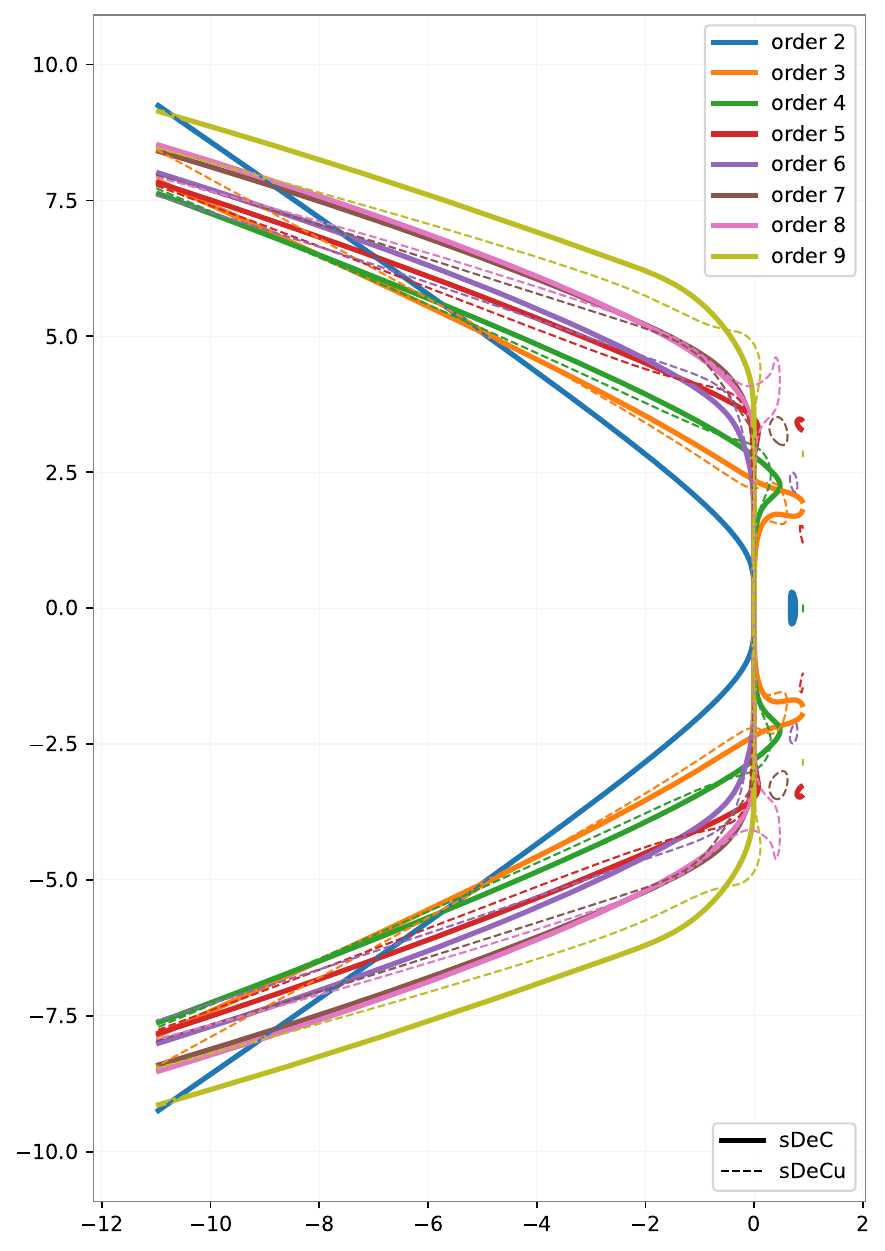}
		\includegraphics[height=0.36\textheight]{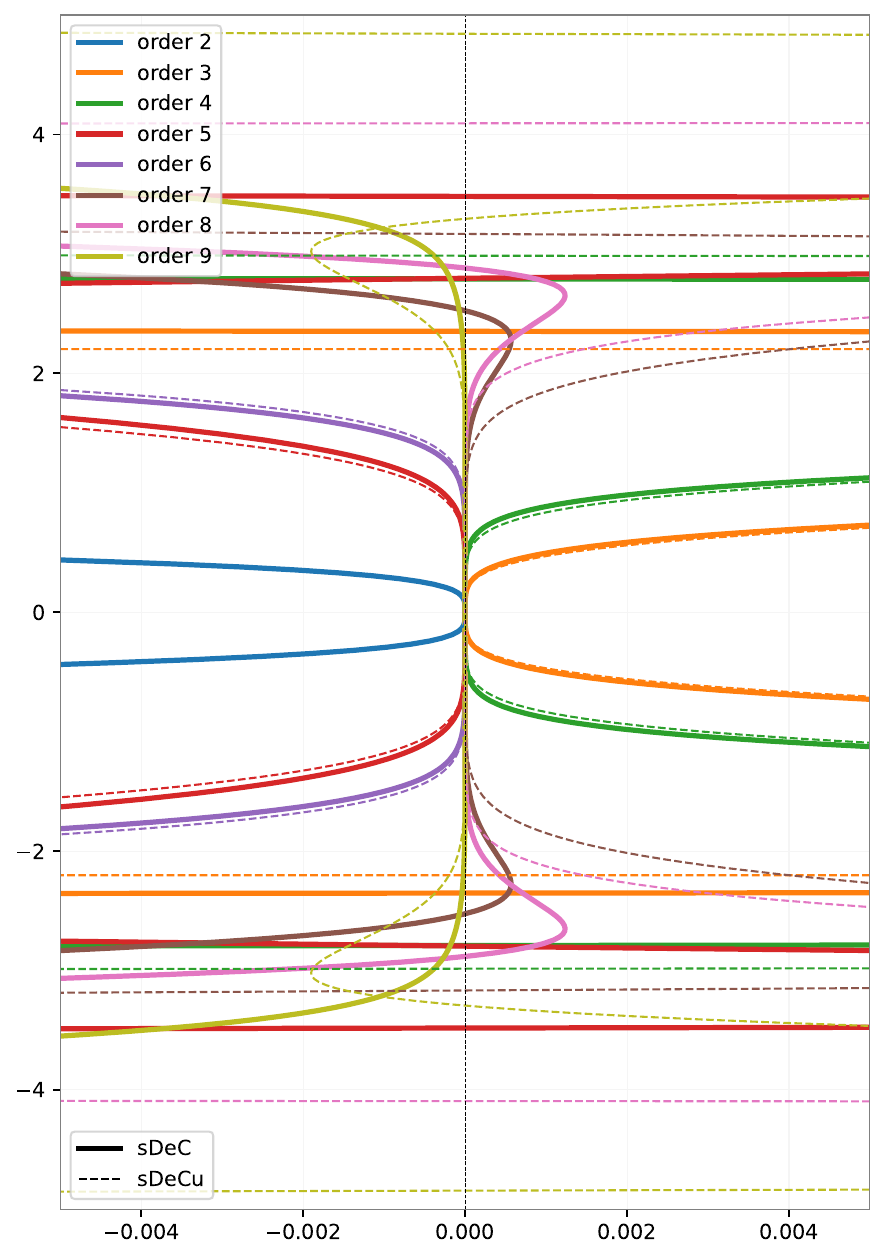}
	\end{center}
	\caption{Minion's stability region for IMEX DeC schemes of orders from 2 to 9 (on the left) with a zoom on the imaginary axis (on the right); $z^{\text{IM}}$ in abscissae, $iz^{\text{EX}}$ in ordinates. We remark that u- and du-approaches coincide on linear problems}
	\label{fig:stability_imex}
\end{figure}

The stability regions of the methods from order 2 to 9 are reported in Figure~\ref{fig:stability_imex} along with some zooms on the imaginary axis. The stability regions are those at the left of the plotted contour lines. 
One can see that the new modified schemes have similar stability regions as the original methods, the only slight difference being observable in the angle of the bDeC methods in the negative real half plane for $z^{\text{IM}}\to -\infty$, with the new methods being slightly less stable. This difference is not present in the sDeC methods.
The zooms on the imaginary axis are of particular interest in the context of applications to pure advection problems with semidiscrete high order methods, whose evolution operators are usually characterized by complex eigenvalues with small real negative part, corresponding to a little amount of (numerical) diffusion. Under this point of view, one can see that both for bDeC and sDeC and related efficient modifications the versions with orders 3, 4, 7, 8 (and 9 only for bDeC)  are more suitable for this kind of problems, as their stability regions effectively contain some portions of the imaginary axis close to the origin with no ``gaps'' with respect to $z^{\text{EX}}$, while the other orders (2, 5 and 6) present some unstable portions of the imaginary axis close to the origin or do not contain at all any part of the imaginary axis, which could be dangerous in the aforementioned context.
Remarkably, in the zooms, bDeCu/du and bDeC coincide. Instead, little differences can be seen between sDeCu/du and sDeC.
In particular, for orders 4, 7 and 8, the efficient modifications seem to bring stability advantages in this context, with sDeCu/du containing a larger portion of the imaginary axis.

\subsection{Fully-implicit stability}\label{sec:stability_fully_implicit}
Here, we consider the stability of the implicit part of the schemes only. Namely, the whole right-hand side is handled implicitly with no explicit terms.
Again, linearity of the schemes implies a generic update of the type $u_{n+1}=R_{\text{FI}}(z)u_{n}$, where $z:=\lambda \dt\in \mathbb C$.
\begin{figure}
	\begin{center}
		Implicit bDeC (continuous), bDeCu/du (dashed)\\
		\includegraphics[height=0.36\textheight]{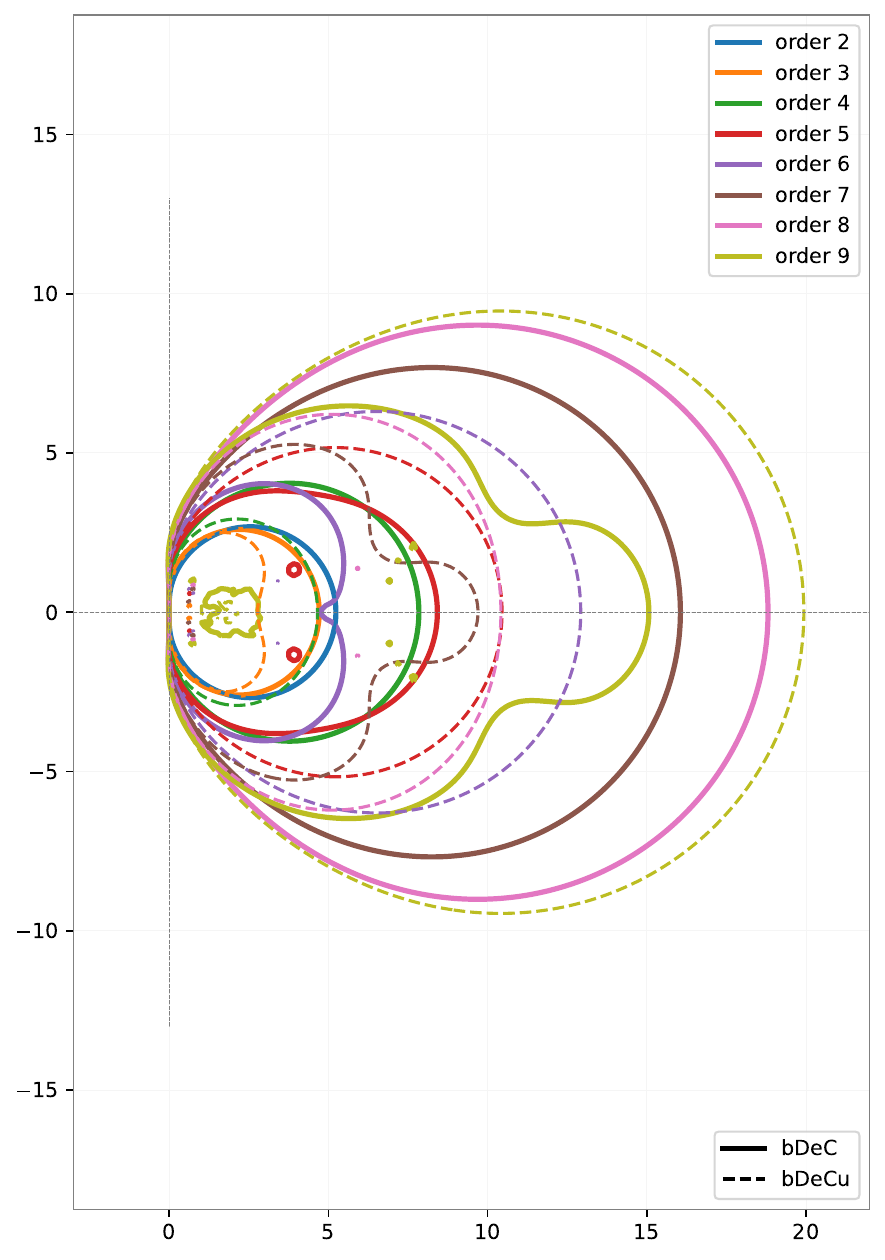}
		\includegraphics[height=0.36\textheight]{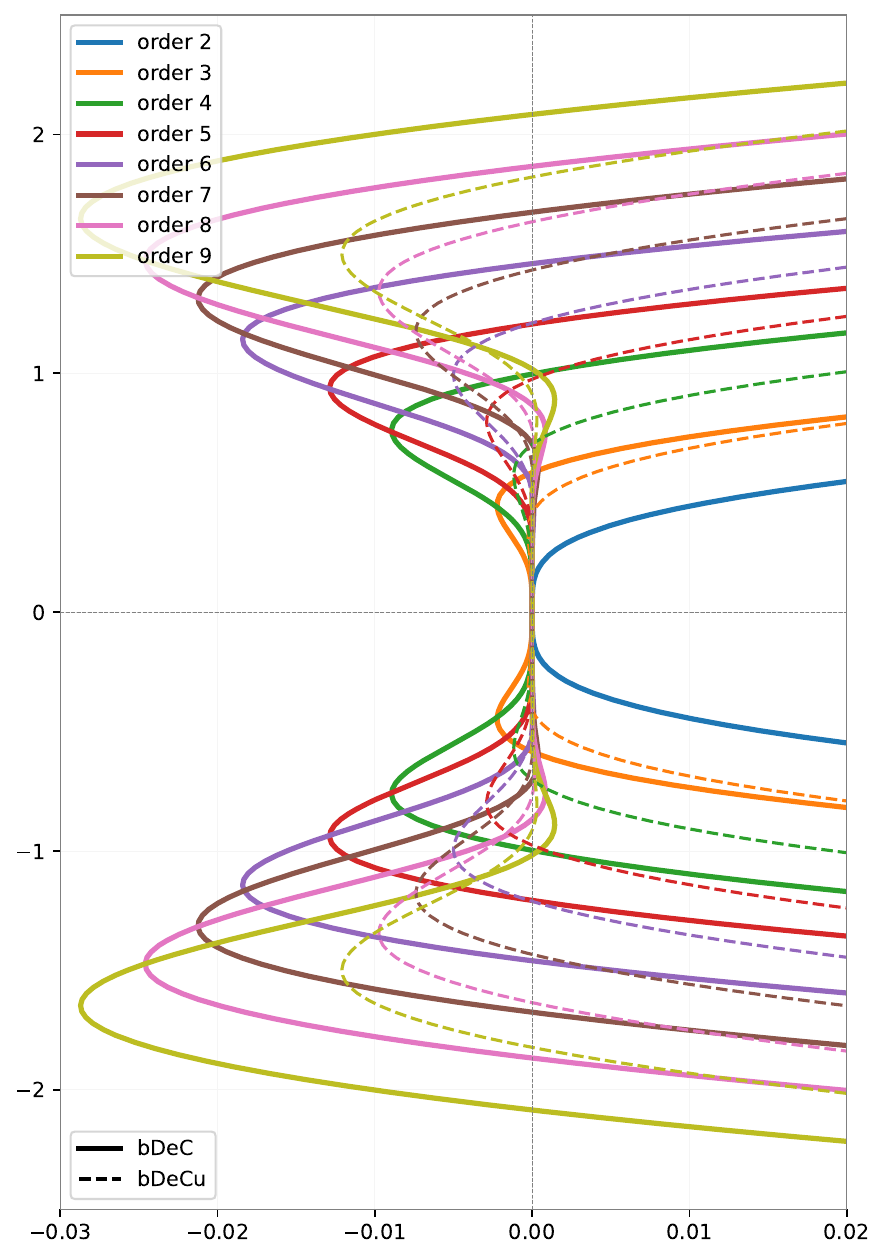}\\ 
		Implicit sDeC (continuous), sDeCu/du (dashed)\\
		\includegraphics[height=0.36\textheight]{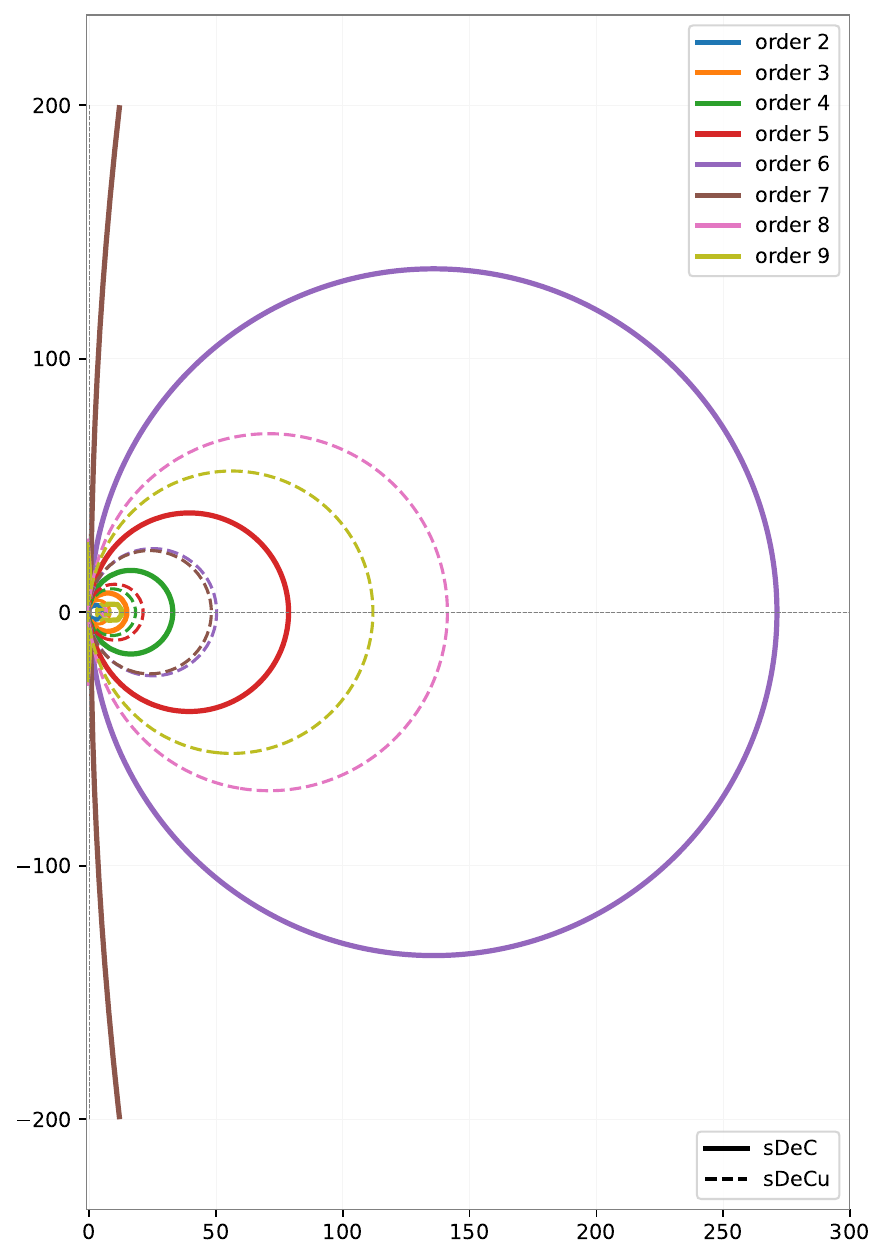}
		\includegraphics[height=0.36\textheight]{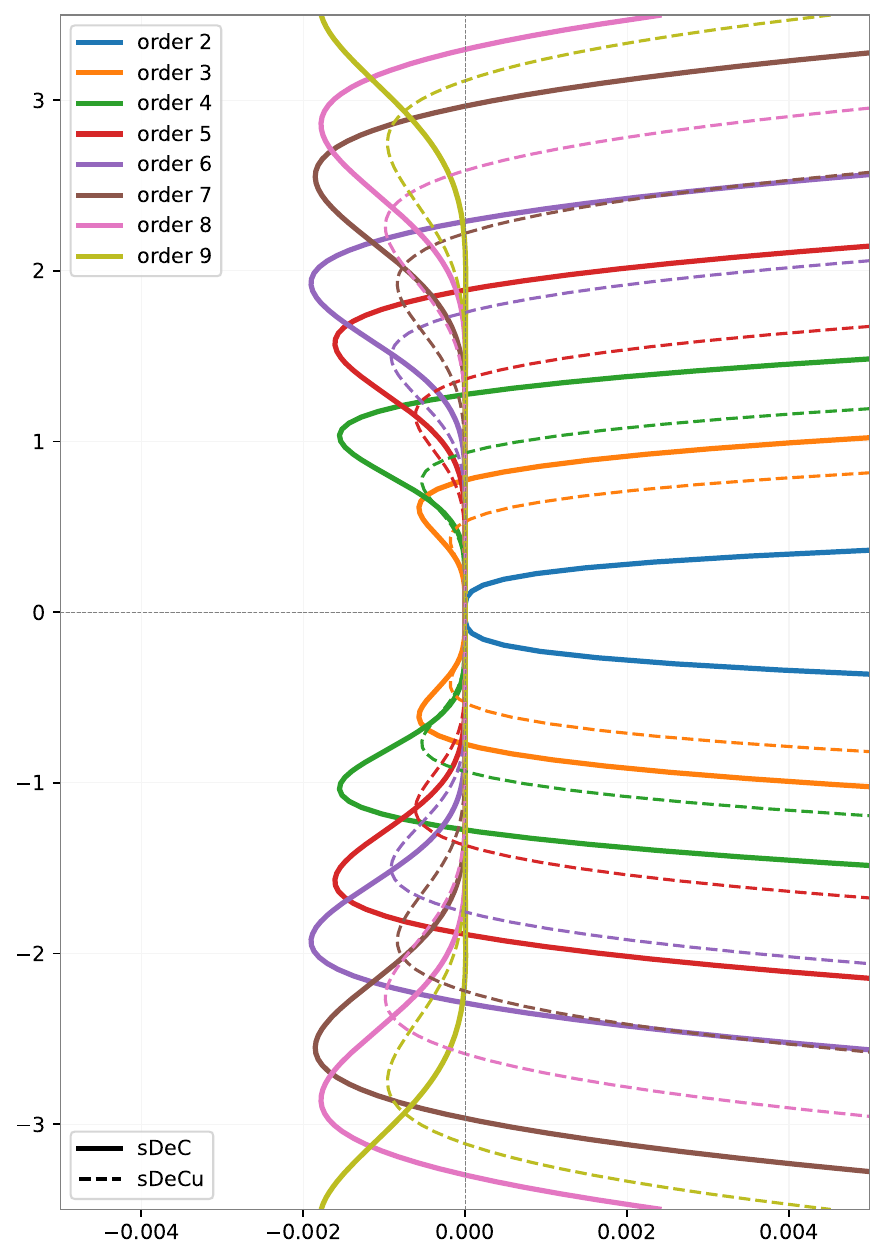}
	\end{center}
	\caption{Stability region for the implicit part of the schemes of orders from 2 to 9 (on the left) with a zoom on the imaginary axis (on the right); $Re(z)$ in abscissae, $Im(z)$ in ordinates. We remark that u- and du-approaches coincide on linear problems}
	\label{fig:stability_implicit}
\end{figure}
In this case, one gets
\begin{equation}\label{eq:fully_implicit_stability_function}
	R_{\text{FI}}(z) = 1+ z\left(\uvec{b}^{\text{IM}}\right)^\top \left(I-zA^{\text{IM}}\right)^{-1} \mathbf{1},
\end{equation}
and the related stability region is
\begin{equation}
	\mathcal{S}_{\text{FI}}:=\left\lbrace z\in \mathbb{C} \quad \text{s.t.} \quad \vert R_{\text{FI}}(z) \vert <1 \right\rbrace.
\end{equation}

The stability regions of the methods from order 2 to 9 are reported in Figure~\ref{fig:stability_implicit} along with some zooms on the imaginary axis.
In this case, one can see that the whole complex half-plane with negative real part is contained in stability regions up to some little localized areas of instability close to the imaginary axis for all orders. 
Remarkably, the introduced interpolation processes reduce the size of such instability regions, and the novel schemes turn out to be more stable than the original versions.

\section{Numerical results}\label{sec:numerical_results}
In this section, we numerically validate the schemes under investigation. To this end, we consider problems of different nature, ranging from smooth tests to assess the order of accuracy, to stiff ODEs and PDE semidiscretizations involving high-order derivatives to assess the ability to handle stiffness.

In Sections~\ref{sec:nonstiff_vibrating_systems} and~\ref{sec:stiff_vibrating_systems}, we consider vibrating systems with moderate and high stiffness, and we test accuracy and computational efficiency of the methods.
In Section~\ref{sec:van_der_pol}, we consider the Van der Pol oscillator problem, to assess the performance of the methods on a nonlinear ODE problem with both stiff and nonstiff regimes.
Sections~\ref{sec:advection_diffusion},~\ref{sec:allen_cahn} and~\ref{sec:cahn_hilliard} are devoted to PDE problems featuring high-order derivatives and nonlinear terms.

\begin{remark}[Nonlinear solver]
	The IMEX DeC formulations require a nonlinear solver to obtain the solution at each subtimenode. In this work, we employ the Newton-Raphson method.
	In principle, one could use the DeC iteration itself as a nonlinear solver (fixed-point iteration), but this would require a number of iterations that is not known a priori and that could be larger than the one required to reach the desired order of accuracy. Experimentally, we have observed that this choice scales very badly for strongly nonlinear problems, while the Newton-Raphson method is more efficient and robust. Hence, for nonlinear problems, we use the Newton-Raphson method to solve every implicit problem with a tolerance of $10^{-12}$ on successive iterations and a maximum of 1000 iterations, which is far away from the average number of iterations usually required (on the order of 10).
\end{remark}
 
\subsection{Moderately stiff vibrating system}\label{sec:nonstiff_vibrating_systems}
Let us consider the following initial value problem
\begin{align}
	\begin{cases}
		m\frac{d^2}{dt^2}y+r\frac{d}{dt}y+ky=F\cos(\Omega t +\varphi),\quad t \in \mathbb{R}^+_0,\\
		y(0)=A ,\\
		\frac{d}{dt}y(0)=B,
	\end{cases}
	\label{eq:nonstiff_vibrating_systems}
\end{align}
with $m,k,\Omega >0$, $r,F,\varphi \geq 0$, corresponding to a damped mechanical vibrating system subject to a sinusoidal external forcing.
The problem can be rewritten as a first order ODE of the type~\eqref{eq:ODE} with the following definitions
\begin{align}
	\uvec{u}=\begin{pmatrix}
		u_1\\
		u_2
	\end{pmatrix}:=\begin{pmatrix}
		y\\
		\frac{d}{dt}y
	\end{pmatrix}, \quad \uvec{N}:=\begin{pmatrix}
		0\\
		-\frac{r}{m}u_2+\frac{F}{m}\cos(\Omega t +\varphi)
	\end{pmatrix}, \quad \uvec{S}:=\begin{pmatrix}
	u_2\\
	-\frac{k}{m}u_1
\end{pmatrix},
\end{align}
where the oscillatory contribution associated with the spring stiffness is treated as a stiff term. Details on how to obtain the exact solution are given in~\cite{micalizzi2023new}.
In this first test, we assume $m=1$, $r=0.01$, $k=10$, $F=1$, $\Omega=2\pi$, $\varphi=\pi/4$, $T_f=10$, $A=0.5$, $B=0.25$.
%
We note that this configuration displays a moderate degree of oscillatory stiffness.


\begin{figure}
	\begin{minipage}{0.49\textwidth}
		\centering
		bDeC
	\end{minipage}
	\begin{minipage}{0.49\textwidth}
		\centering sDeC
	\end{minipage}\\
	\includegraphics[width=0.49\textwidth]{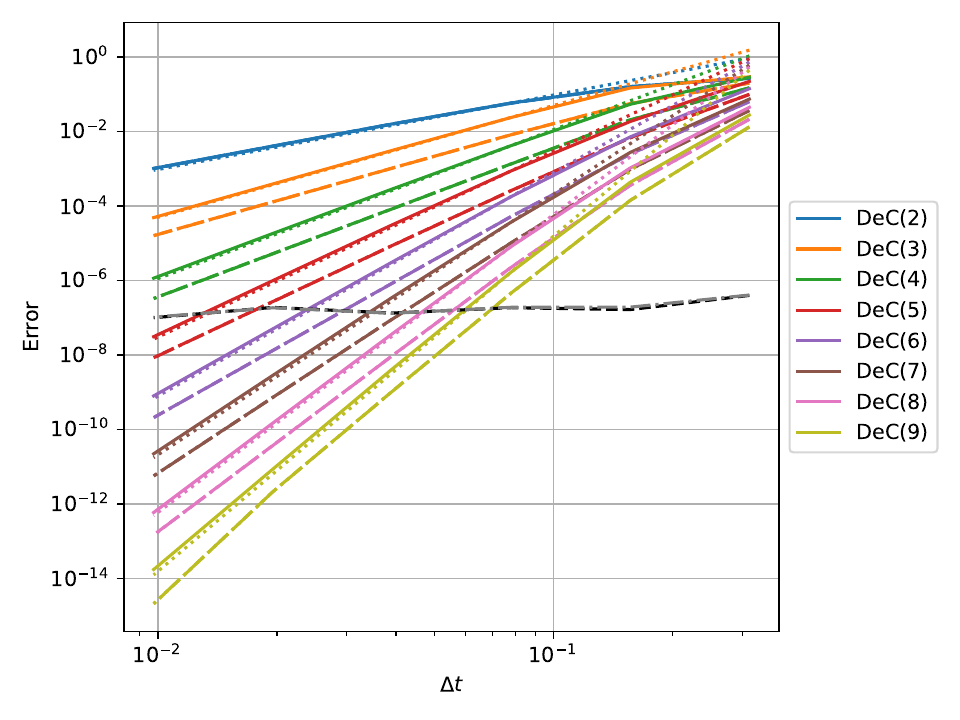}
	\includegraphics[width=0.49\textwidth]{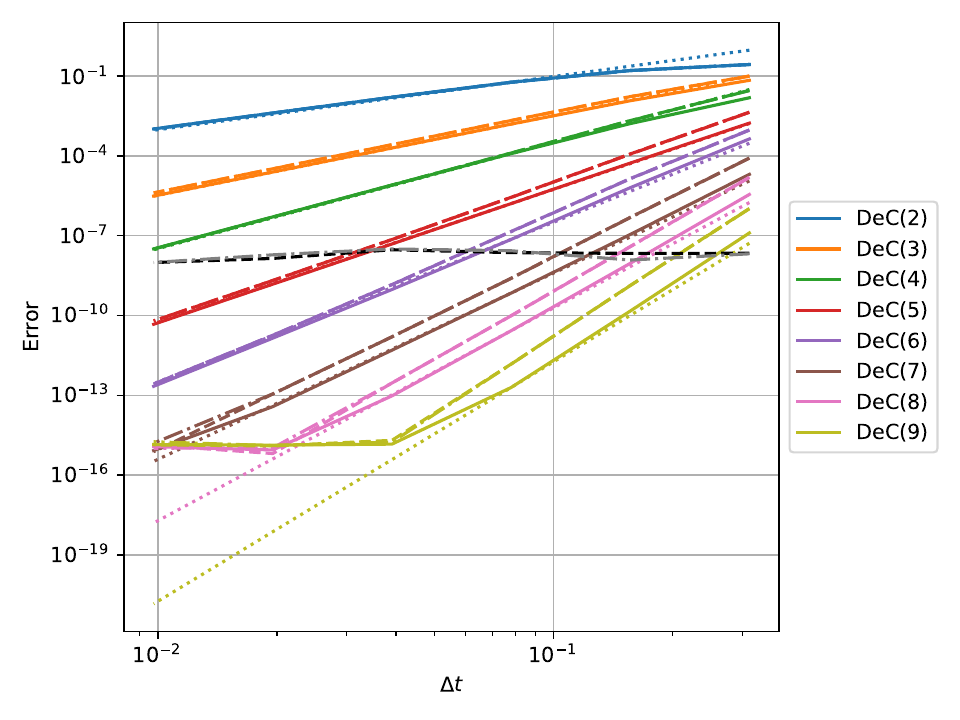}\\
	\includegraphics[width=0.49\textwidth]{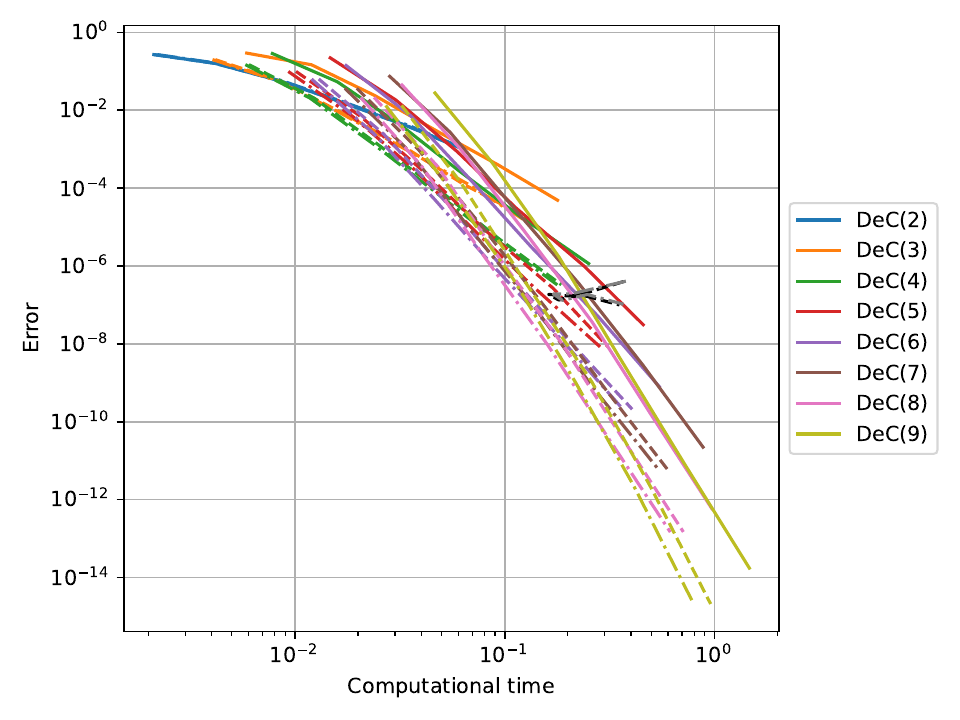}
	\includegraphics[width=0.49\textwidth]{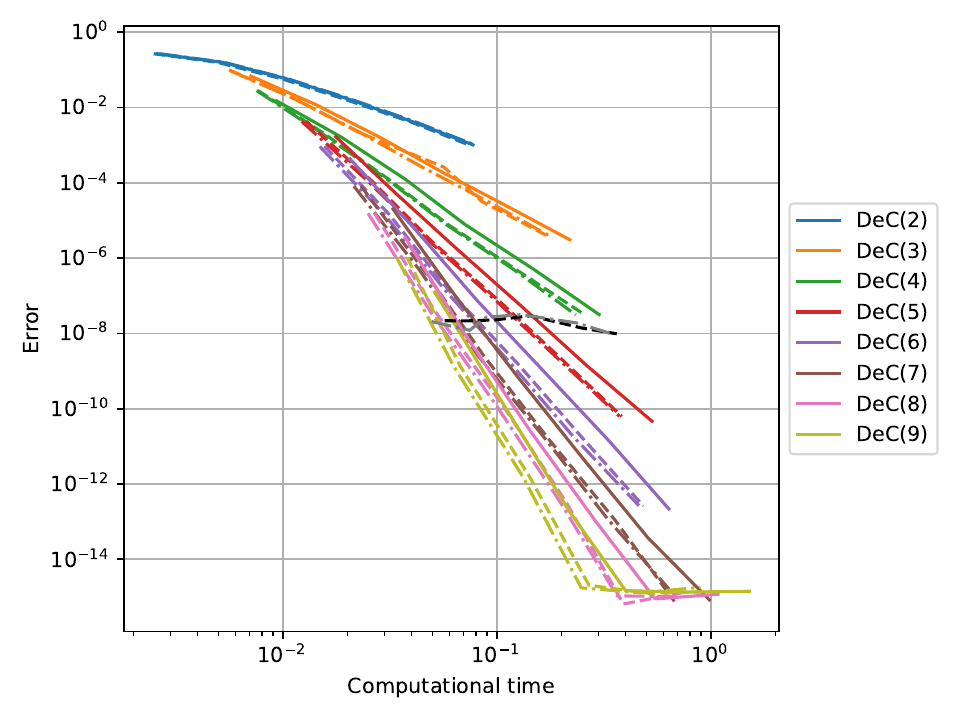}\\
	\caption{Moderately stiff vibrating system: Convergence and efficiency analysis of the IMEX DeC schemes. Top: convergence; Bottom: efficiency. Left: bDeC; Right: sDeC. Original DeC methods are depicted in continuous lines, the DeCu in dashed lines and the DeCdu in dashed-dotted lines. The dotted lines in the first two plots represent the expected order of accuracy. The black and gray horizontal lines represent the adaptive methods obtained with a tolerance of $\varepsilon=10^{-7}$, black for DeCu, gray for DeCdu. }
	\label{fig:convergence_vibrating}
\end{figure}

In Figure~\ref{fig:convergence_vibrating}, we report the errors of the methods against time-step size and computational time. The number of time-steps is taken between 32 and 1024. The results show that the expected order of accuracy is achieved for all methods. Furthermore, bDeCu and bDeCdu are indistinguishable in the convergence plots as well as sDeCu and sDeCdu. 
Noticeably, while sDeCu and sDeCdu produce higher errors for fixed $\dt$ with respect to the original sDeC, in the bDeC case the novel schemes produce smaller errors with respect to the original version. In all cases, however, one can appreciate how the novel modifications are computationally more efficient than the original schemes, in particular for high order.
As one can see and as noticed also in \cite{offner2025analysis}, the bDeC schemes in the implicit form are way less accurate than the sDeC schemes, in particular, in the high order case; this also applies to their efficient modifications.
On the other hand, we remark that bDeC schemes, contrarily to sDeC ones, can be parallelized and be extremely competitive~\cite{ketcheson2014comparison}.

\begin{figure}
	\centering
\includegraphics[width=0.8\textwidth]{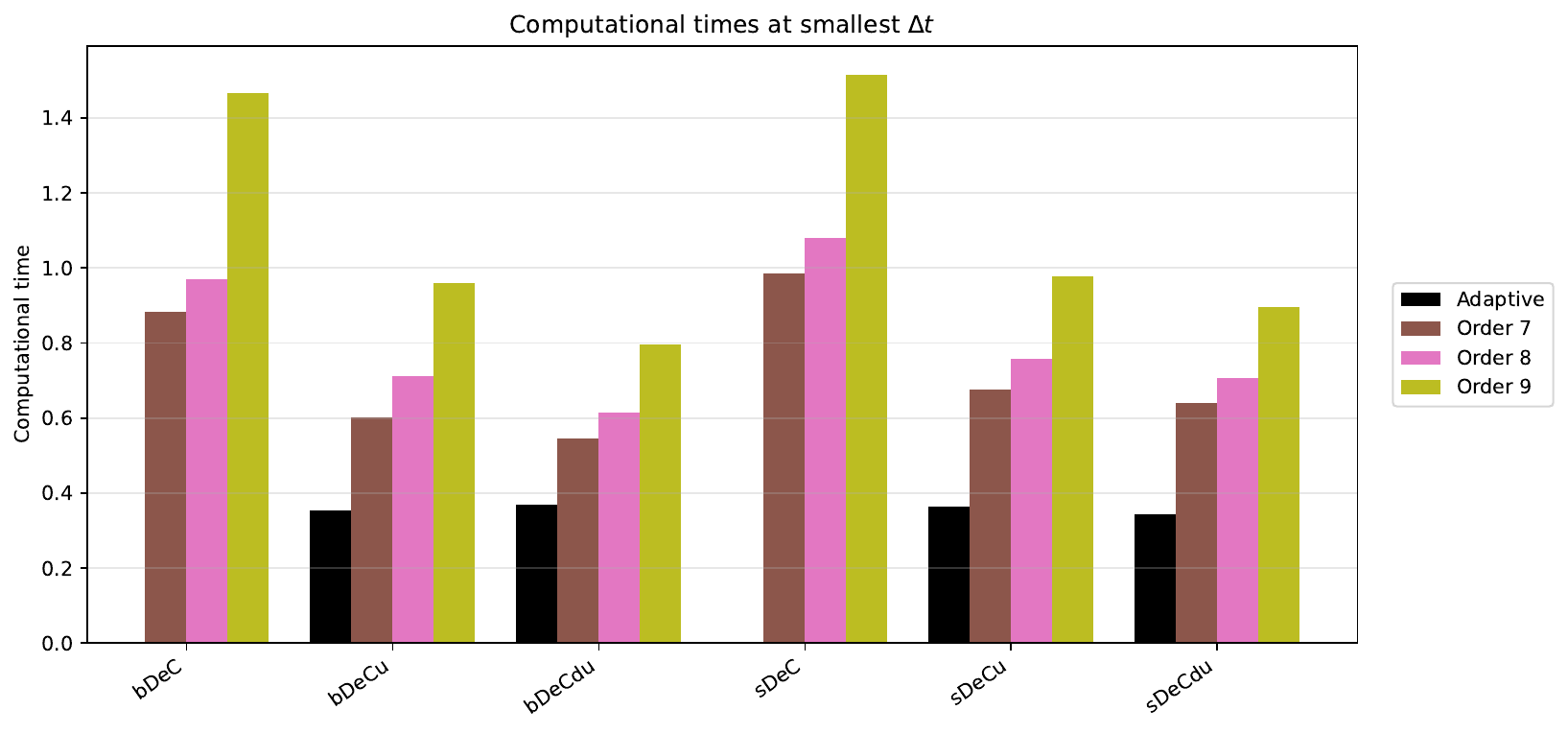}
\caption{Moderately stiff vibrating system: Computational times corresponding to the smallest $\dt$ for the adaptive methods and the non-adaptive schemes of orders 7, 8, and 9.}
\label{fig:barplot_time_vibrating}
\end{figure}

\begin{figure}
	\centering
	\includegraphics[width=0.8\textwidth]{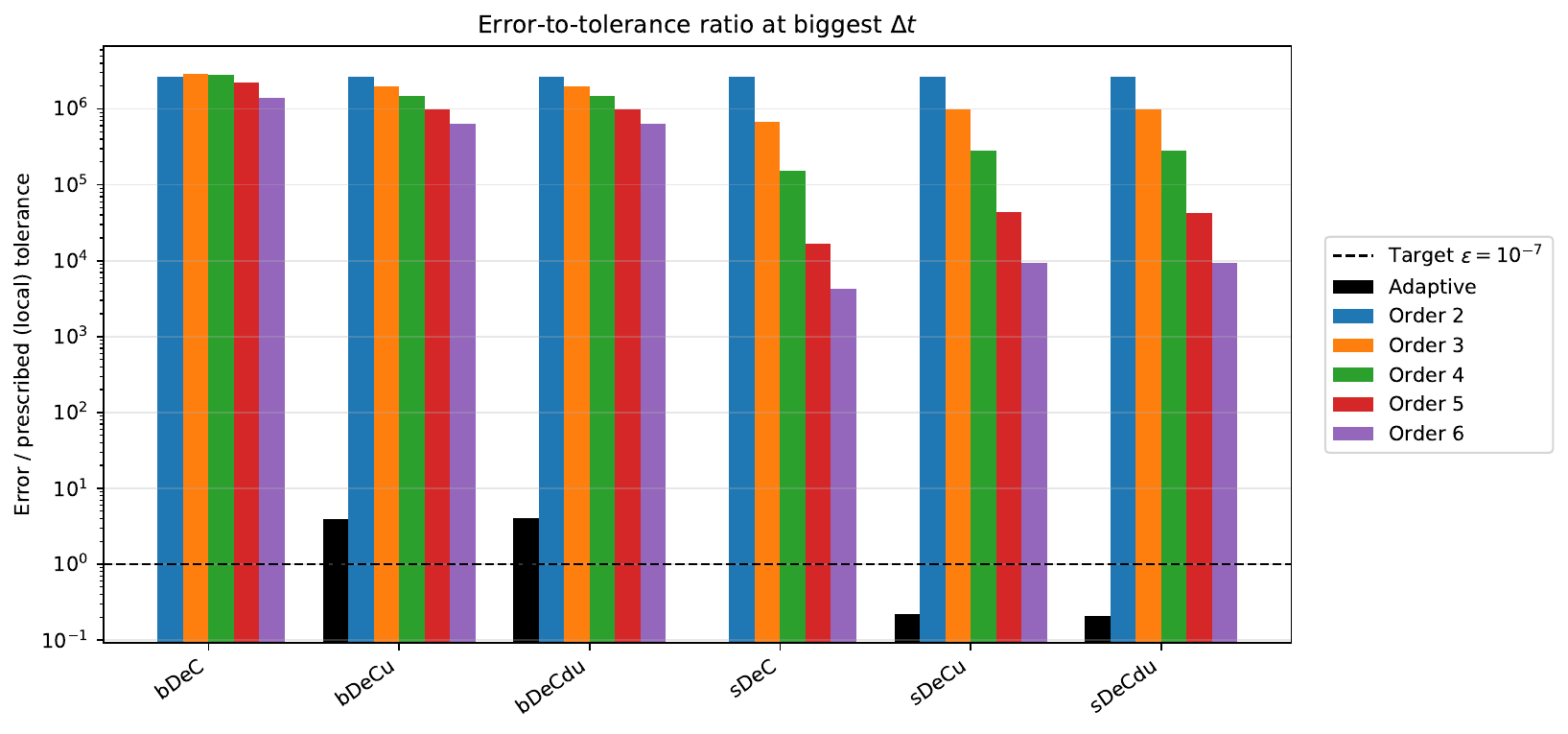}
	\caption{Moderately stiff vibrating system: Ratio between the error and the prescribed (local) tolerance corresponding to the biggest $\dt$ for the adaptive methods and the non-adaptive schemes of orders from 2 to 6.}
	\label{fig:barplot_error_vibrating}
\end{figure}

In gray and black, we depict the two adaptive versions obtained, according to the strategy described in Section~\ref{sec:adaptivity}, with a tolerance of $\varepsilon=10^{-7}$. The adaptive methods are able to reach an error consistent with the prescribed accuracy independently of the time-step size.
To further highlight the advantages of the adaptive strategy, in Figure~\ref{fig:barplot_time_vibrating} we report the computational times corresponding to the smallest $\dt$ for the adaptive methods and for all the considered non-adaptive schemes of orders 7, 8, and 9.
The adaptive methods effectively adjust the order of accuracy according to the prescribed tolerance, leading to a substantial reduction in computational cost. In contrast, fixed very-high-order schemes are characterized by a considerably larger computational cost, although such a high order is not required to meet the prescribed accuracy tolerance.
In Figure~\ref{fig:barplot_error_vibrating}, instead, we report the ratios corresponding to the biggest $\dt$ between the error and the tolerance used for the adaptive simulations for the adaptive methods and all the considered non-adaptive schemes from order 2 to 6.
As one can see, the adaptive strategy is able to achieve a final error of the same order of magnitude as the prescribed tolerance, while fixed low-order schemes are characterized by a much higher error.
Let us notice that the adaptation strategy is local and this makes the final error slightly higher than the prescribed (local) tolerance. The issue can be solved by prescribing a stricter tolerance through a safety coefficient and/or through estimates of the final error from the local one.
In any case, the results indicate that the adaptation strategy is sufficiently robust across a wide range of refinements.

\subsection{Stiff vibrating system}\label{sec:stiff_vibrating_systems}
\begin{figure}
	\centering
	\includegraphics[width=0.8\textwidth]{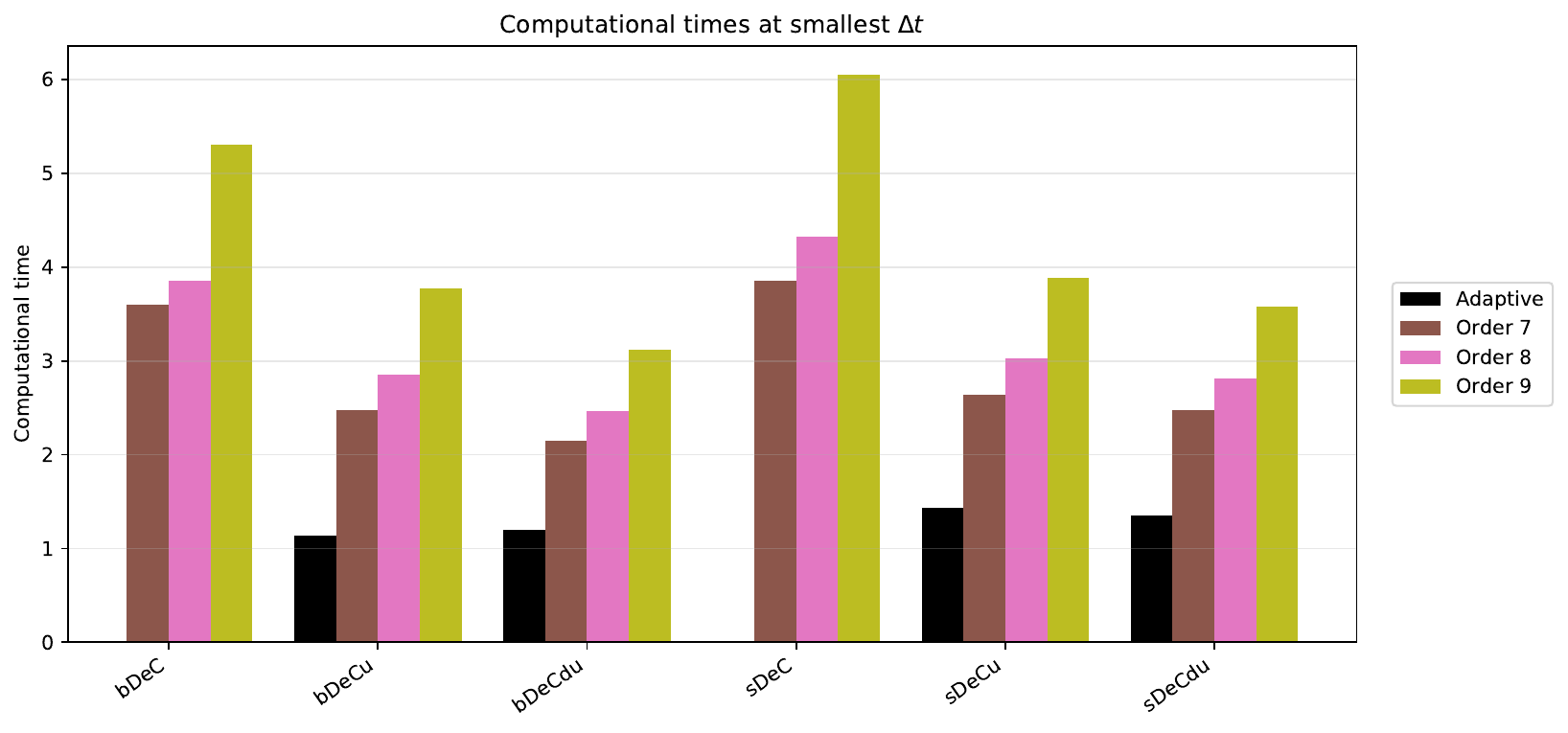}
	\caption{Stiff vibrating system: Computational times corresponding to the smallest $\dt$ for the adaptive methods and the non-adaptive schemes of orders 7, 8, and 9.}
	\label{fig:barplot_time_vibrating_high}
\end{figure}

\begin{figure}
	\centering
	\includegraphics[width=0.8\textwidth]{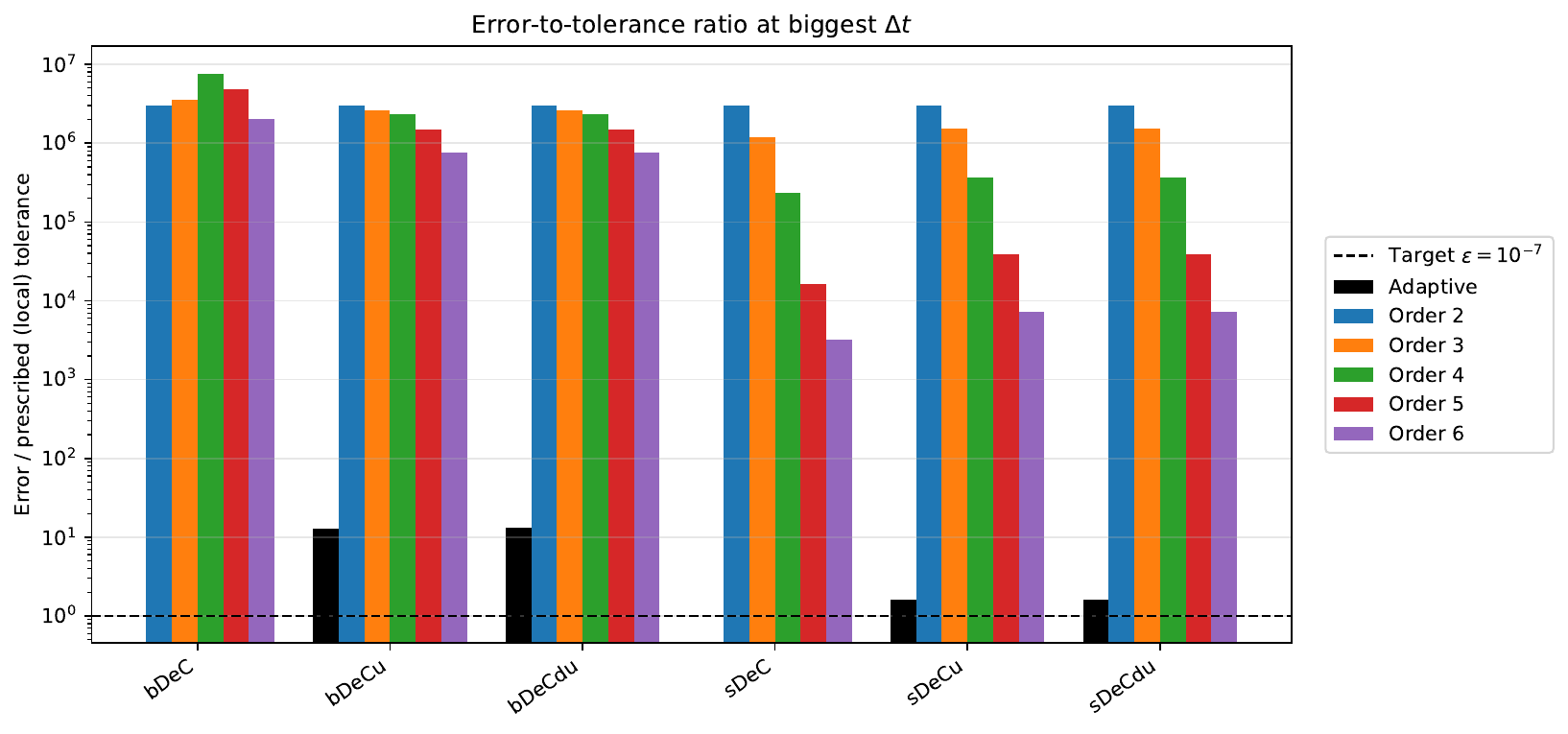}
	\caption{Stiff vibrating system: Ratio between the error and the prescribed (local) tolerance corresponding to the biggest $\dt$ for the adaptive methods and the non-adaptive schemes of orders from 2 to 6.}
	\label{fig:barplot_error_vibrating_high}
\end{figure}

Now, we consider the same problem as in the last section but with a much higher level of stiffness with $k=100$.
\begin{figure}
	\centering
	\begin{minipage}{0.49\textwidth}
		\centering
		bDeC
	\end{minipage}
	\begin{minipage}{0.49\textwidth}
		\centering sDeC
	\end{minipage}\\
	\includegraphics[width=0.49\textwidth]{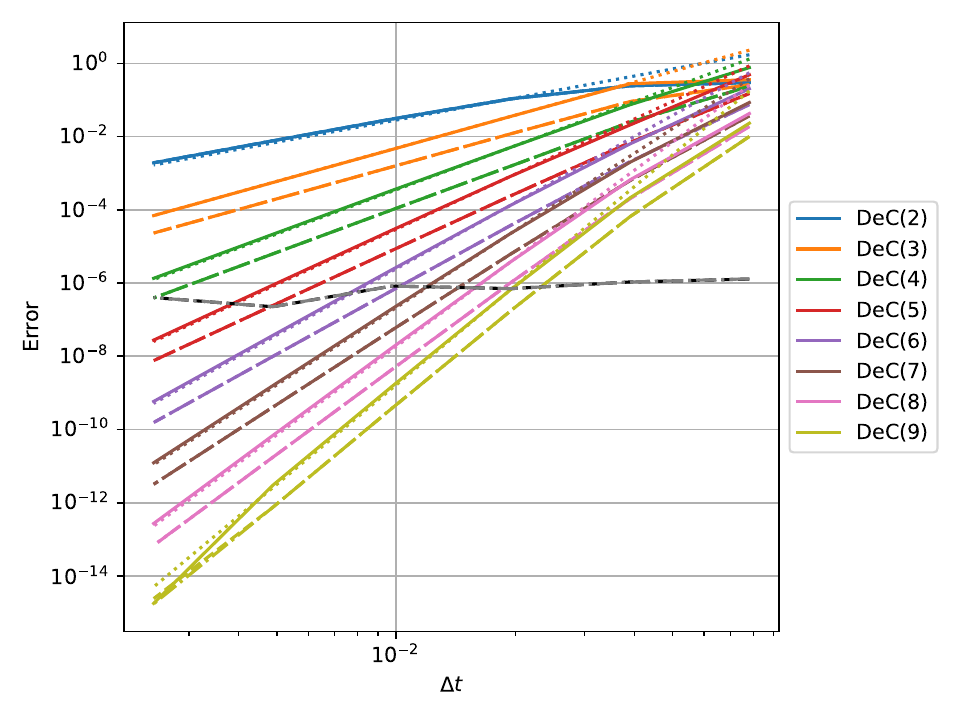}
	\includegraphics[width=0.49\textwidth]{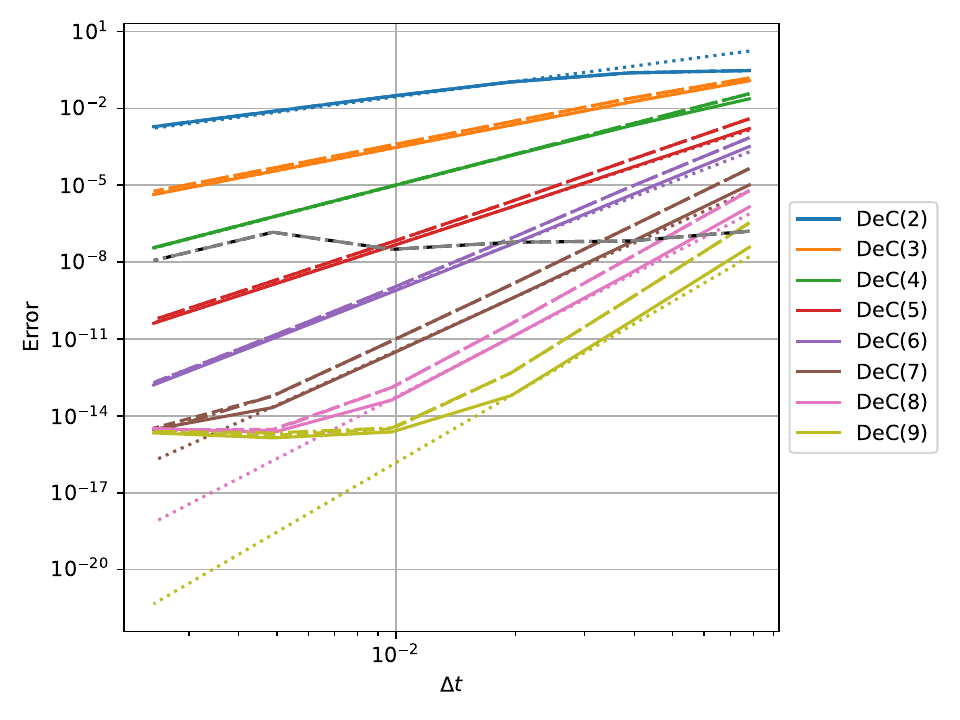}\\
	\includegraphics[width=0.49\textwidth]{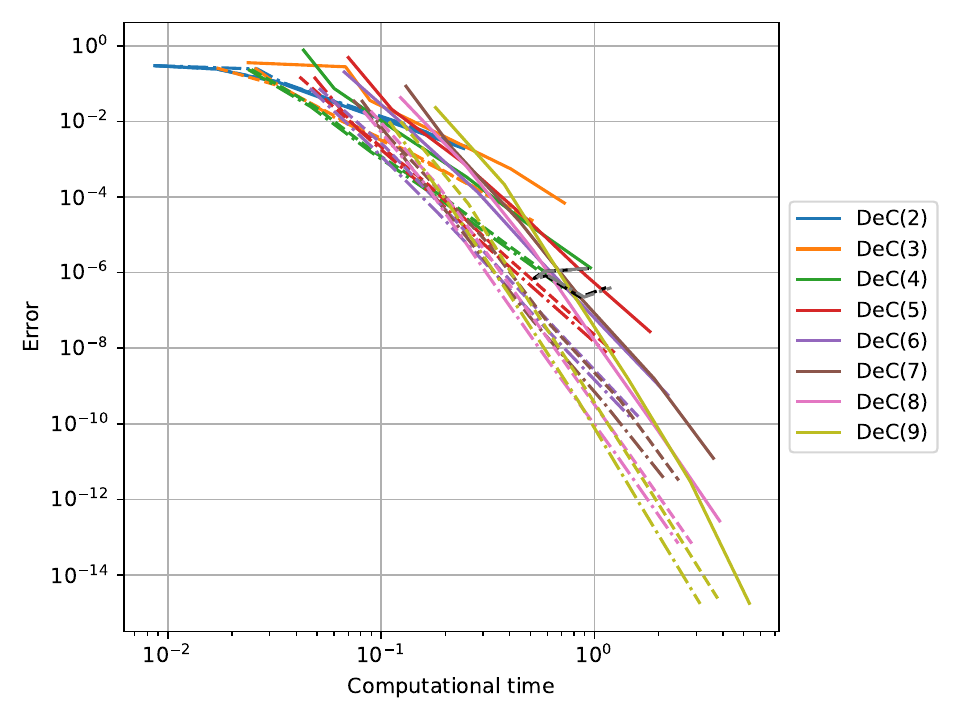}
	\includegraphics[width=0.49\textwidth]{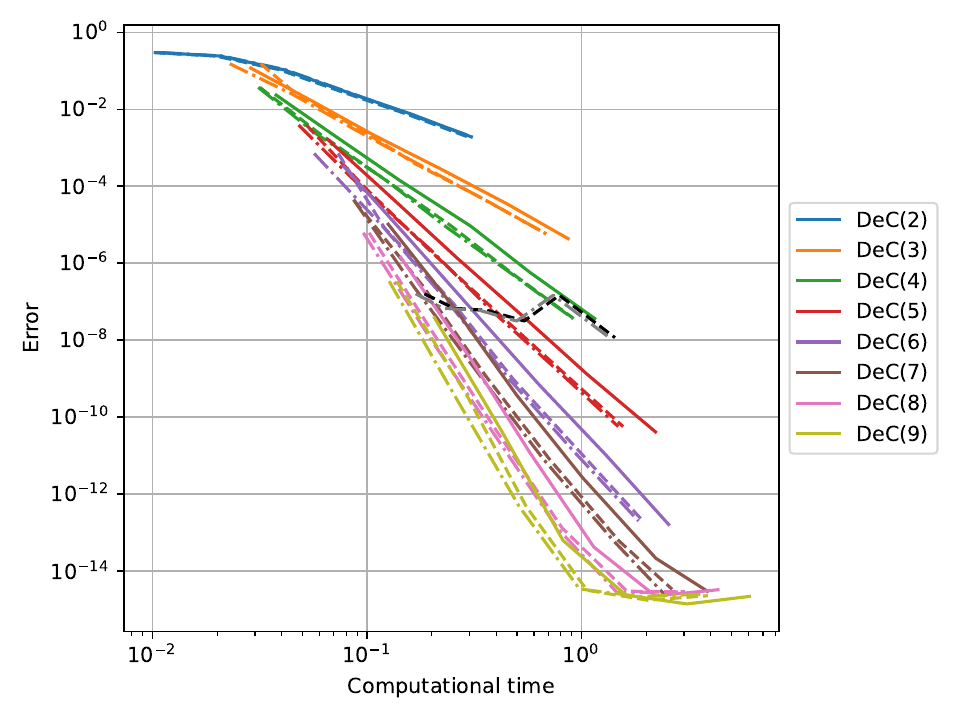}\\
	\caption{Stiff vibrating system: Convergence and efficiency analysis of the IMEX DeC schemes. Top: convergence; Bottom: efficiency. Left: bDeC; Right: sDeC. Original DeC methods are depicted in continuous lines, the DeCu in dashed lines and the DeCdu in dashed-dotted lines. The dotted lines in the first two plots represent the expected order of accuracy. The black and gray horizontal lines represent the adaptive methods obtained with a tolerance of $\varepsilon=10^{-7}$, black for DeCu, gray for DeCdu.}
	\label{fig:convergence_vibrating_high}
\end{figure}
To run the convergence analysis for this test, we refined the time mesh, starting from 128 time-steps up to 4096 time-steps.  
In Figure~\ref{fig:convergence_vibrating_high}, we report the errors of the methods against time-step size and computational time. The results show that the expected order of accuracy is achieved for all methods. 
Same considerations as for the previous test apply.
For fixed $\dt$, sDeCu and sDeCdu produce bigger errors with respect to sDeC, while bDeCu and bDeCdu produce smaller errors with respect to bDeC.
In all cases, the modified methods are computationally convenient with respect to the original ones, especially for high order.
Finally, the adaptive schemes, still with $\varepsilon=10^{-7}$, are able to adaptively select the order of accuracy to match the prescribed tolerance.
Again, the advantages of employing the adaptive strategy can be appreciated by comparing the computational costs at smallest $\dt$ in Figure~\ref{fig:barplot_time_vibrating_high} and the error-to-tolerance ratios at biggest $\dt$ in Figure~\ref{fig:barplot_error_vibrating_high}.

\subsection{Van der Pol oscillator}\label{sec:van_der_pol}
In this test, we consider the Van der Pol oscillator, introduced in~\cite{van1926lxxxviii} to model nonlinear oscillations in electrical circuits, rescaled in the form presented in \cite{wanner1996solving,boscarino2018implicit}
\begin{align}
	\begin{cases}
		& \frac{d^2}{dt^2}y=\mu \left[(1-y^2)\frac{d}{dt}y-y\right],\\
		& y(0)=A ,\\
		& \frac{d}{dt}y(0)=B,
	\end{cases}
	\label{eq:van_der_pol}
\end{align}
where $\mu>0$ is a constant.
Also in this case, we can rewrite the problem in form~\eqref{eq:ODE} with 
\begin{align}
	\uvec{u}=\begin{pmatrix}
		u_1\\
		u_2
	\end{pmatrix}:=\begin{pmatrix}
		y\\
		\frac{d}{dt}y
	\end{pmatrix}, \quad \uvec{N}:=\uvec{0}, \quad \uvec{S}:=\begin{pmatrix}
		u_2\\
		\mu \left[(1-u_1^2)u_2-u_1\right]
	\end{pmatrix},
\end{align}
corresponding to a fully-implicit treatment of the right-hand side, as it is not trivial to split the system into a fast and a slow dynamic part. 
We set $\mu=5$, $A=2$, $B=-\frac{2}{3}+\frac{10}{810}\frac{1}{\mu^2}-\frac{292}{2187}\frac{1}{\mu^4}$ and $T_f=20$.

\begin{figure}
	\begin{minipage}{0.435\textwidth}
	\centering
	Order 5 $N=30$\\
	\includegraphics[width=\textwidth, trim={0 0 120 0},clip]{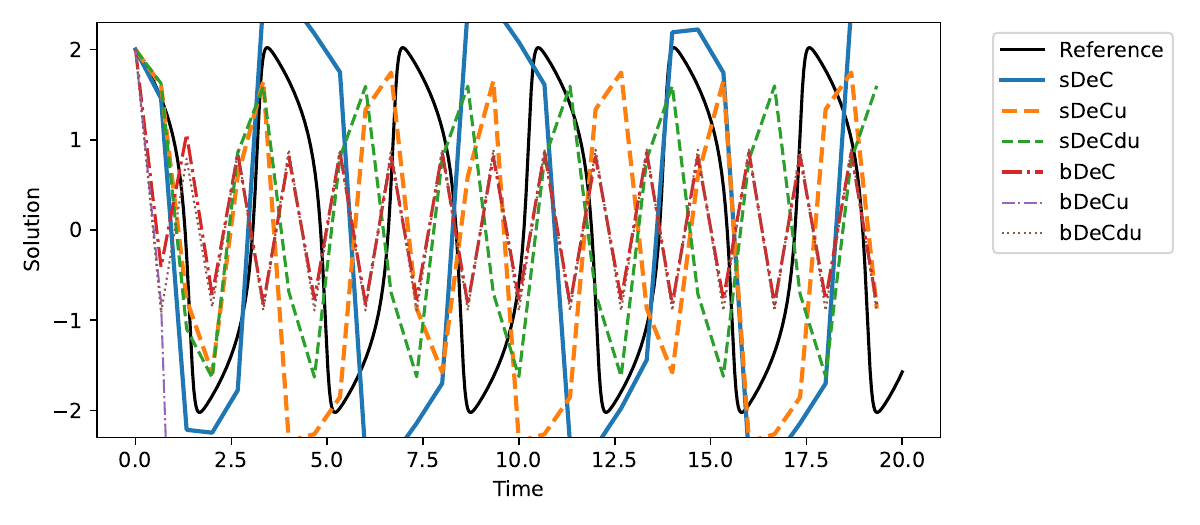}\\
	Order 5 $N=150$\\
	\includegraphics[width=\textwidth, trim={0 0 120 0},clip]{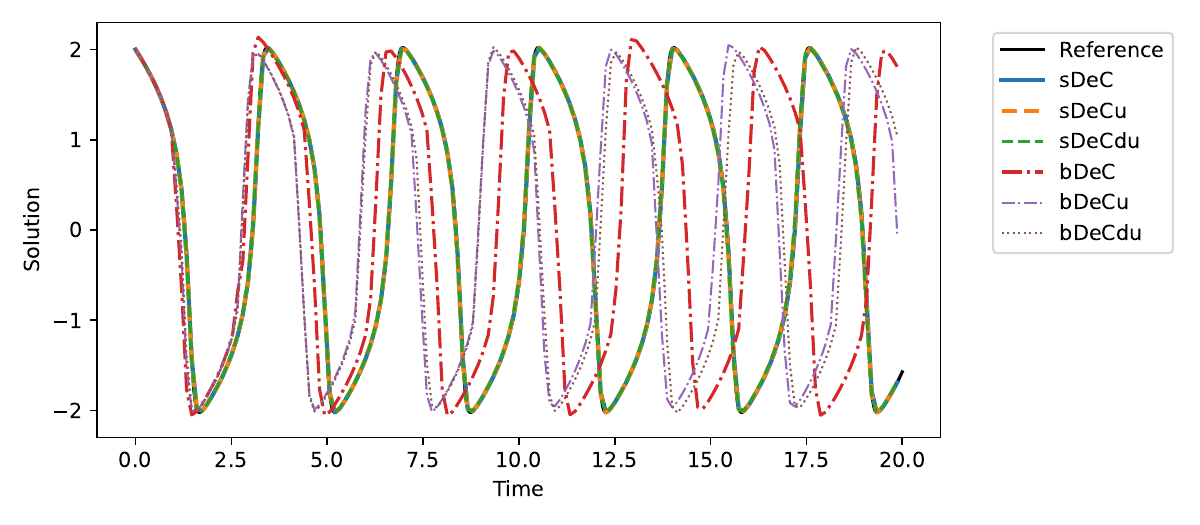}		
	\end{minipage}
	\begin{minipage}{0.435\textwidth}
	\centering
	Order 5 $N=80$\\
	\includegraphics[width=\textwidth, trim={0 0 120 0},clip]{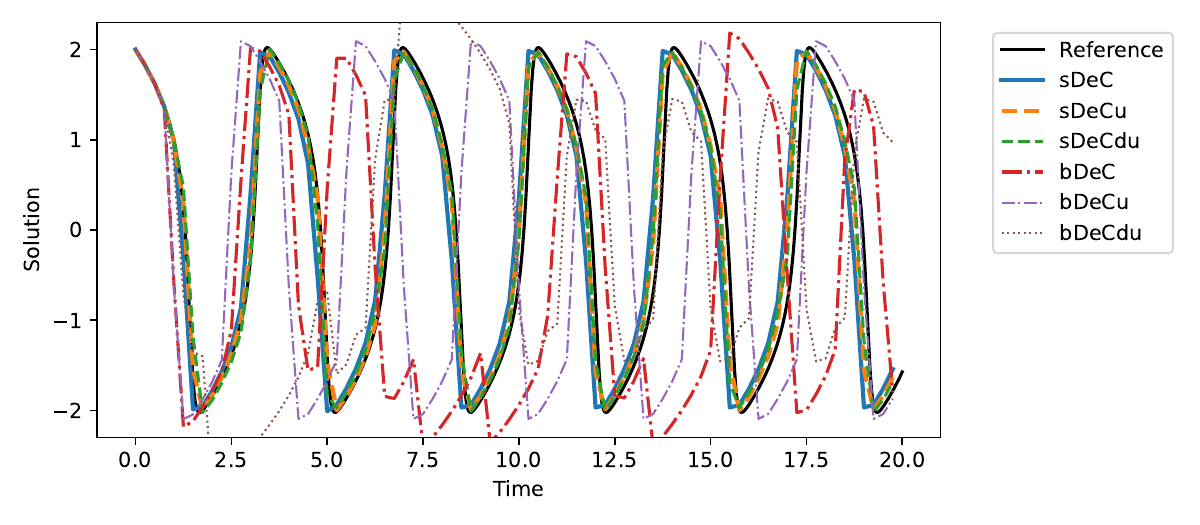}\\
	Order 5 $N=300$\\
	\includegraphics[width=\textwidth, trim={0 0 120 0},clip]{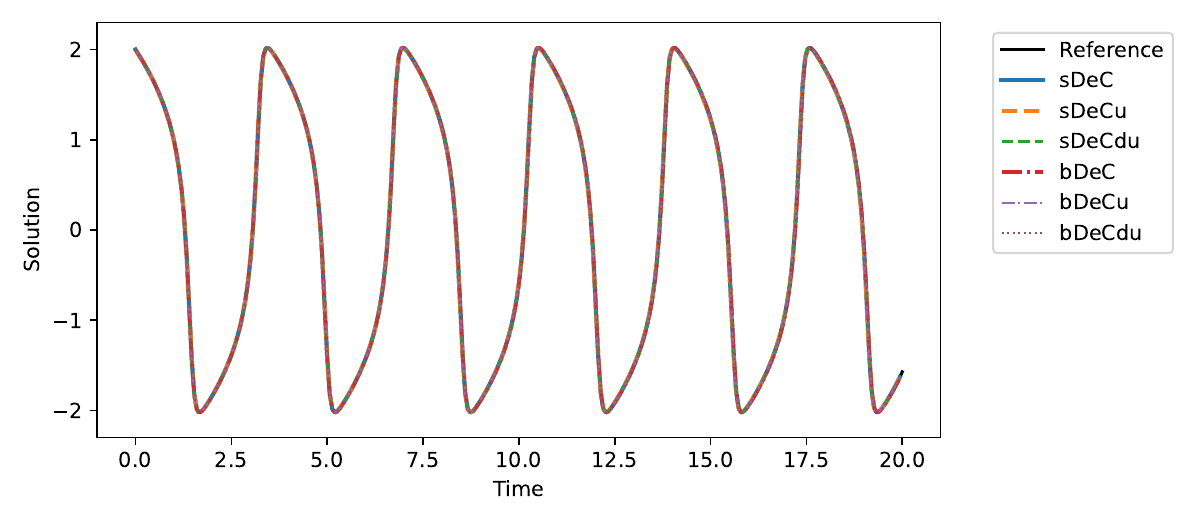}\\
	\end{minipage}
	\begin{minipage}{0.11\textwidth}
	\centering
	\includegraphics[width=\textwidth, trim={470 0 10 0},clip]{RescaledVanDerPol_solution_order_5_N_80.pdf}
	\end{minipage}
	\caption{Van der Pol oscillator: Results obtained with fifth-order methods for different temporal resolutions} \label{fig:vdp_order5}
\end{figure}
\begin{figure}
	\begin{minipage}{0.435\textwidth}
	\centering
	$N=30$\\
	\includegraphics[width=\textwidth, trim={0 0 140 30},clip]{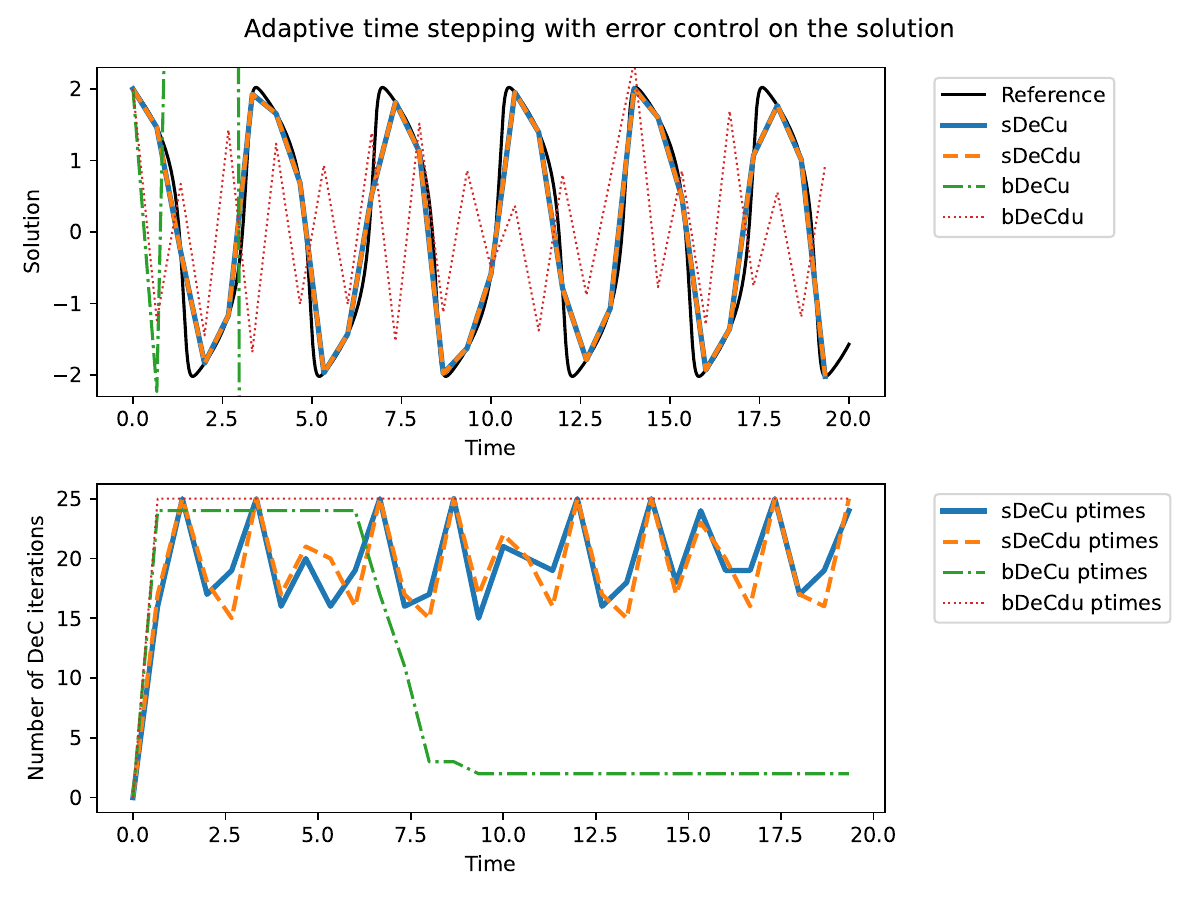}\\
	$N=150$\\
	\includegraphics[width=\textwidth, trim={0 0 140 30},clip]{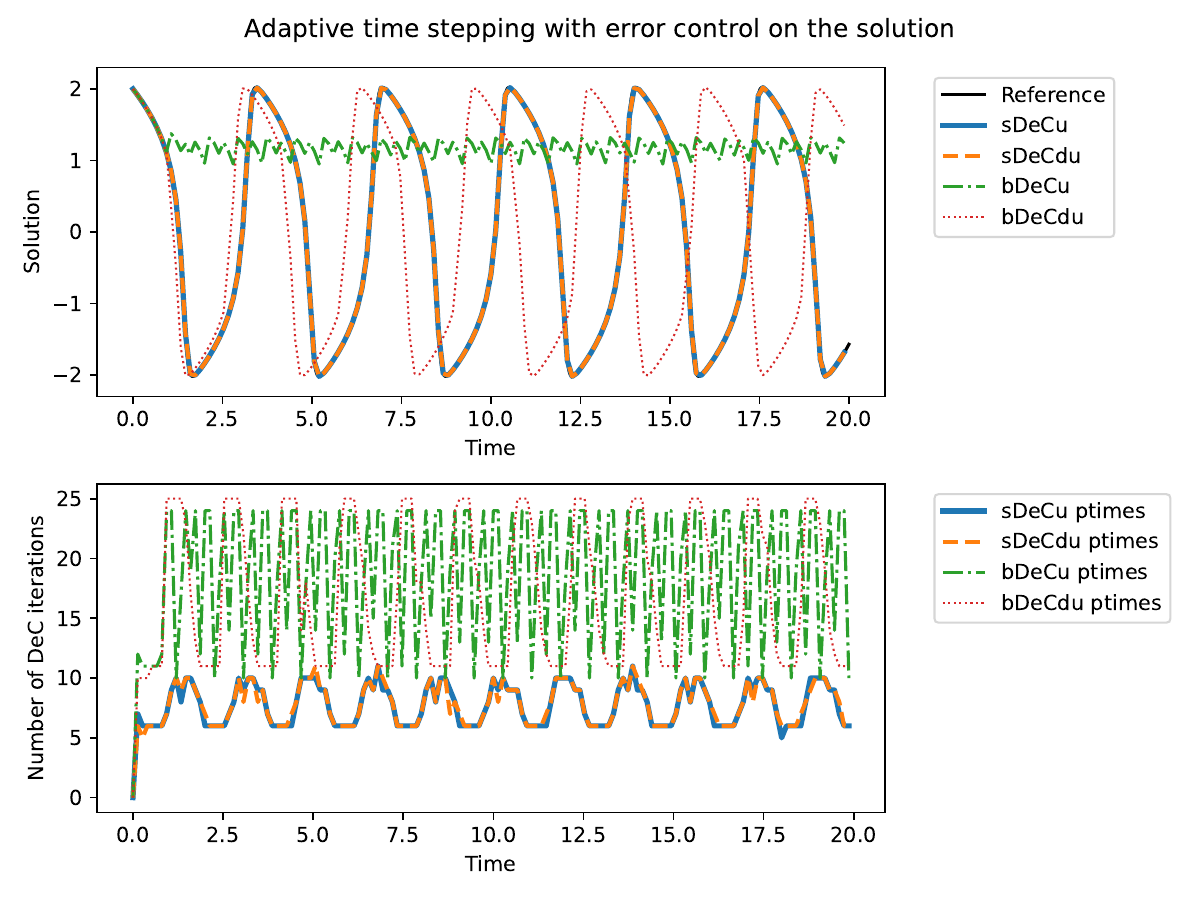}		
	\end{minipage}
	\begin{minipage}{0.435\textwidth}
	\centering
	$N=80$\\
	\includegraphics[width=\textwidth, trim={0 0 140 30},clip]{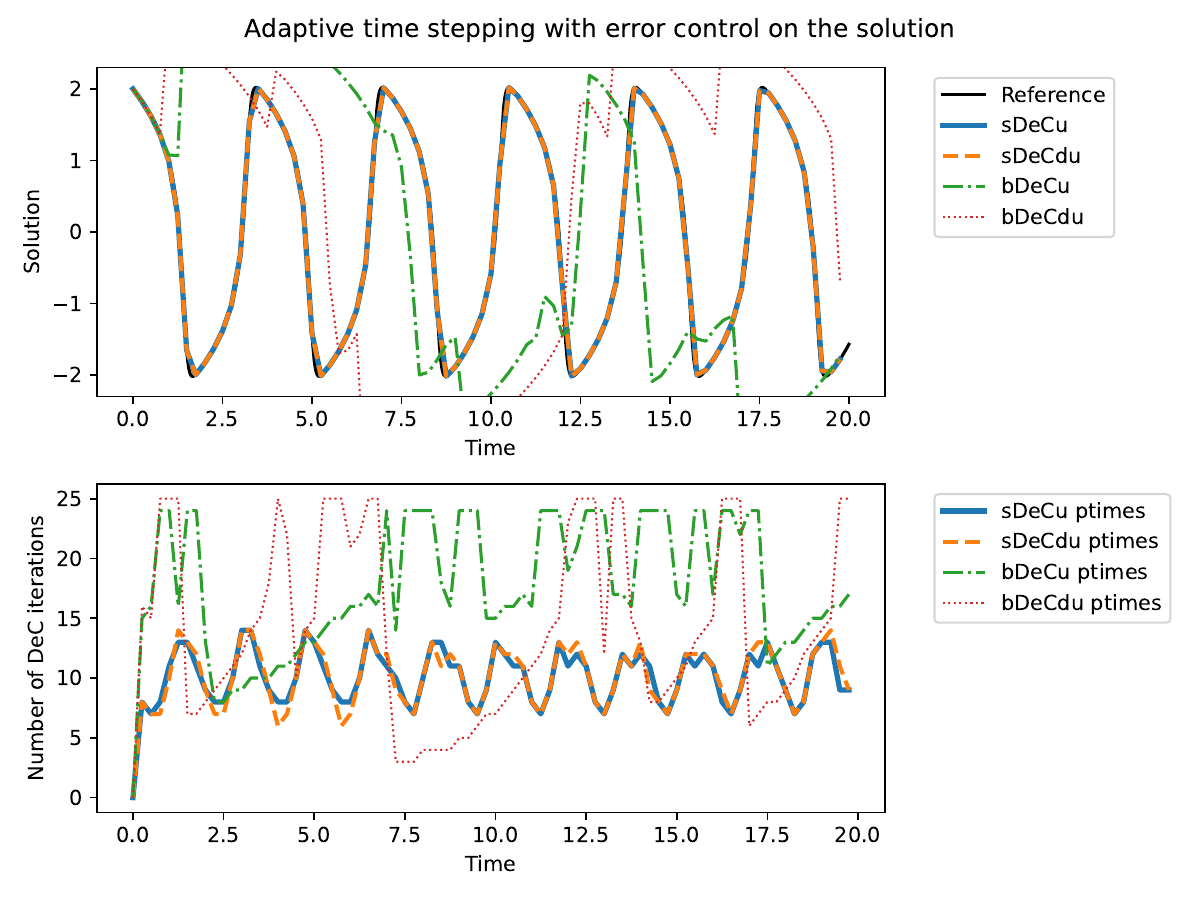}\\
	$N=300$\\
	\includegraphics[width=\textwidth, trim={0 0 140 30},clip]{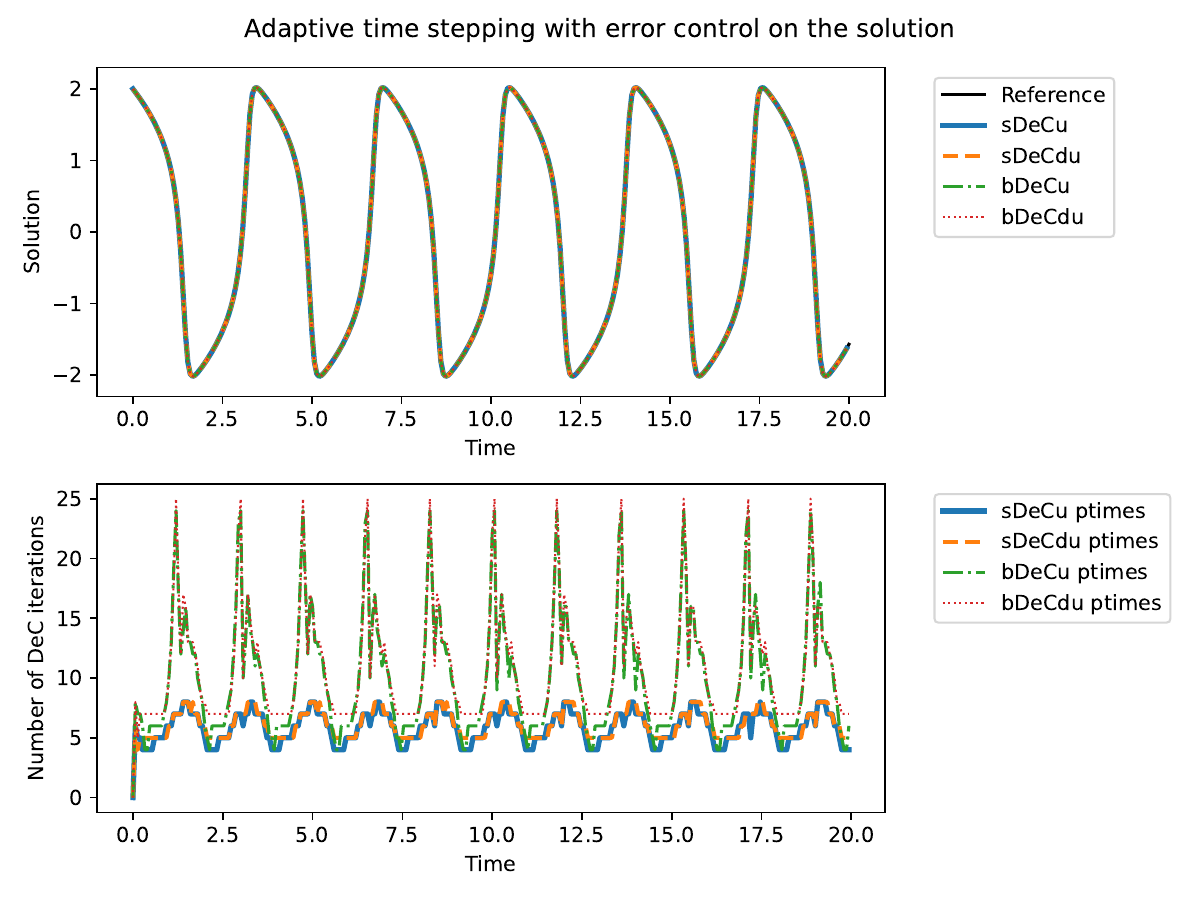}\\
	\end{minipage}
	\begin{minipage}[t]{0.11\textwidth}
	\centering
	\includegraphics[width=\textwidth, trim={440 200 30 20},clip]{RescaledVanDerPol_solution_adaptive_N_80.pdf}\\
	\end{minipage}
	\caption{Van der Pol oscillator: Results obtained with adaptive methods with tolerance $\varepsilon=10^{-6}$ and maximum number of iterations allowed equal to 25 for different temporal resolutions. Top: numerical solutions; bottom: number of DeC iterations used in each time-step}\label{fig:vdp_adapt}
\end{figure}
In Figure~\ref{fig:vdp_order5}, we report the solution obtained with the fifth-order methods for different time-step sizes. It is clear that, even with an implicit treatment, it is necessary to use a sufficiently fine temporal resolution to capture the dynamics of the stiff problem. In particular, for $N=30$ all methods have difficulties in correctly simulating the dynamics. For $N=80$ the sDeC methods are able to capture the solution, while the bDeC methods still fail to do so. Even at $N=150$, the bDeC methods do not accurately reproduce the oscillation frequency, while for $N=300$ all the schemes capture the solution.

In Figure~\ref{fig:vdp_adapt}, we report the results obtained with the adaptive methods. In this context, we fix a maximum order of accuracy equal to 25, hence fixing also the maximum number of subtimenodes to $14$, and perform at most 25 iterations if the convergence condition on the last subtimenode in Equation~\eqref{eq:tol} is not matched. The number of subtimenodes is increased throughout the iterative procedure until $14$, and the last iterations are performed with such a configuration. The tolerance is set to $\varepsilon:=10^{-6}$.

Again, the sDeC methods show better performance than the bDeC methods. In particular, even for $N=30$, they are able to correctly capture the solution, using between 15 and 25 iterations at each time-step, choosing more iterations when the problem is stiffer.
On the other hand, the bDeC methods struggle even for finer temporal resolutions. What we have observed is that the Newton solver reaches a solution without ambiguity, but the DeC (Picard) iteration process often oscillates between two different solutions (in the last subtimenode), hence not satisfying the stopping criterion. 
This indicates that, for this test and parameter regime, the bDeC methods are less robust than the sDeC methods.
Again, we remark that a suitable parallelization of bDeC schemes and related adaptive versions might allow running simulations for smaller $\dt$ with similar efficiency performance, but this aspect is not investigated here.

\subsection{Advection--diffusion equation}\label{sec:advection_diffusion}
Now, we move to PDE tests. We start from the one--dimensional advection--diffusion equation
\begin{equation}\label{eq:advection_diffusion}
	u_t+au_x=\mu u_{xx}, \quad (x,t)\in \Omega\times[0,T_f],
\end{equation}
with $a\in\mathbb{R}$ and $\mu\geq 0$ being the advection and diffusion coefficients respectively, and $\Omega:=[x_L,x_R]$ and $T_f>0$ being the spatial domain and the final time.
To discretize the spatial operators, we employ a finite-difference framework with arbitrary order of accuracy.

We consider a uniform spatial grid $\Omega_{\Delta x}=\left\{x_j \ : \ x_j = x_L + j\Delta x, \ j \in \{0,\hdots , J\}, \ \Delta x:=\frac{x_R-x_L}{J} \right\}$, and we denote by $u_j$ the approximation of $u(x_j,t)$.

To approximate the advection term in \eqref{eq:advection_diffusion}, namely the first spatial derivative $\partial_x u(x)$, we use the stable finite-difference stencils introduced in \cite{Iserles1982}. Thus, $\partial_x u$ is approximated at $x_j$ by an $[r,s]$-discretization,
\begin{equation}\label{eq:r-s_scheme}
\partial^{[r,s]}_{\Delta x}(u(x_j)) = \frac{1}{\Delta x} \sum\limits_{k=-r}^s \alpha_k u_{j+k},
\end{equation}
with $r,s$ chosen such that $\alpha_{-r}, \alpha_{s} \neq 0$.
The highest attainable order of such a discretization is $q=r+s$, and the coefficients in \eqref{eq:r-s_scheme} are uniquely determined by imposing that the discretization is of order $q$:
\begin{align}\label{eq:def_advection_optimal_stencil}
\alpha_0&=
\begin{cases}
-\sum\limits_{k=r+1}^s\frac{1}{k}, & s\ge r+1, \\
0, & s=r, \\
\sum\limits_{k=s+1}^r\frac{1}{k}, & r\ge s+1,
\end{cases} \qquad
\alpha_k = \frac{(-1)^{k+1}}{k}\cdot \frac{r!s!}{(r+k)!(s-k)!}, \quad -r\le k \le s, \ k\neq0.
\end{align}
As shown in \cite{Iserles1982}, these so-called optimal-order schemes of order $q$ are stable if and only if $s \le r \le s+2$ for $a>0$.
We make use of these stable optimal-order schemes in our analysis. In particular, as we assume $a>0$, for order $q\geq 1$ we adopt the following upwinded stencil choice: $r=\left\lceil \frac{q+1}{2}\right\rceil$ and $s:=q-r$.
%
%
For the diffusion term, we use a central finite-difference approximation of the second spatial derivative $\partial_{xx}u(x)$, as reported in Table~\ref{tab:CFD_schemes_second_deriv}.
In particular, for each order, we use the smallest stencil able to guarantee the desired accuracy. 
After spatial discretization, the advection contribution is treated
explicitly, while the diffusion contribution is treated implicitly.

\begin{table}
	\centering
	\caption{Central finite-difference approximations of $\partial_{xx}u$ at the grid point $x_j$ \cite{fornberg_finite_difference}}
	\label{tab:CFD_schemes_second_deriv}
	\scriptsize
	\resizebox{0.8\textwidth}{!}{%
		\begin{tabular}{|c|c|}
			\hline
			order & finite-difference approximation of $\partial_{xx}u(x_j)$\\
			\hline
			2 &
			$\displaystyle
			\frac{1}{\Delta x^2}
			\left(
			u_{j-1}-2u_j+u_{j+1}
			\right)$
			\\
			\hline
			4 &
			$\displaystyle
			\frac{1}{\Delta x^2}
			\left(
			-\frac{1}{12}u_{j-2}
			+\frac{4}{3}u_{j-1}
			-\frac{5}{2}u_j
			+\frac{4}{3}u_{j+1}
			-\frac{1}{12}u_{j+2}
			\right)$
			\\
			\hline
			6 &
			$\displaystyle
			\frac{1}{\Delta x^2}
			\left(
			\frac{1}{90}u_{j-3}
			-\frac{3}{20}u_{j-2}
			+\frac{3}{2}u_{j-1}
			-\frac{49}{18}u_j
			+\frac{3}{2}u_{j+1}
			-\frac{3}{20}u_{j+2}
			+\frac{1}{90}u_{j+3}
			\right)$
			\\
			\hline
			8 &
			$\displaystyle
			\frac{1}{\Delta x^2}
			\left(
			-\frac{1}{560}u_{j-4}
			+\frac{8}{315}u_{j-3}
			-\frac{1}{5}u_{j-2}
			+\frac{8}{5}u_{j-1}
			-\frac{205}{72}u_j
			+\frac{8}{5}u_{j+1}
			-\frac{1}{5}u_{j+2}
			+\frac{8}{315}u_{j+3}
			-\frac{1}{560}u_{j+4}
			\right)$
			\\
			\hline
		\end{tabular}
	}
\end{table}

\begin{figure}
	\centering
	Order 5\\
	\includegraphics[width=0.49\textwidth]{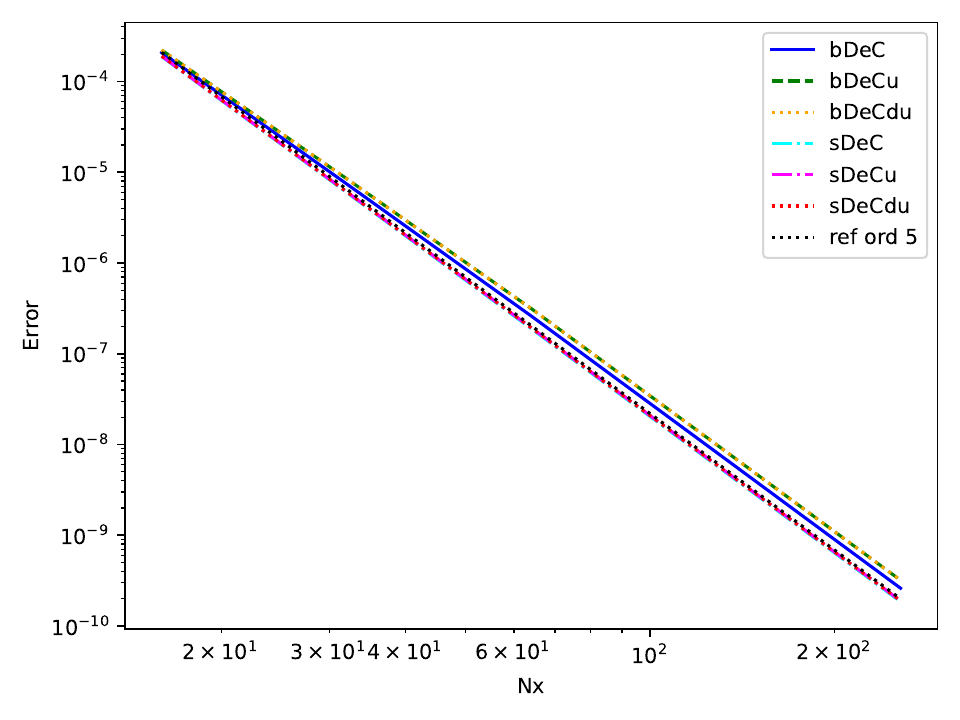}
	\includegraphics[width=0.49\textwidth]{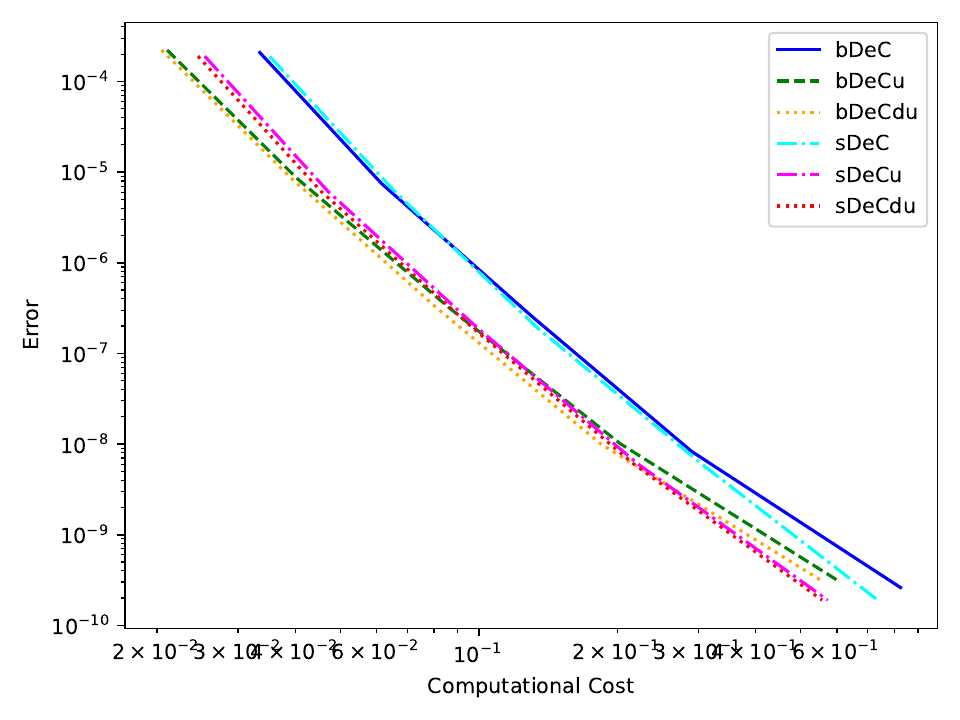}\\
	Order 8\\
	\includegraphics[width=0.49\textwidth]{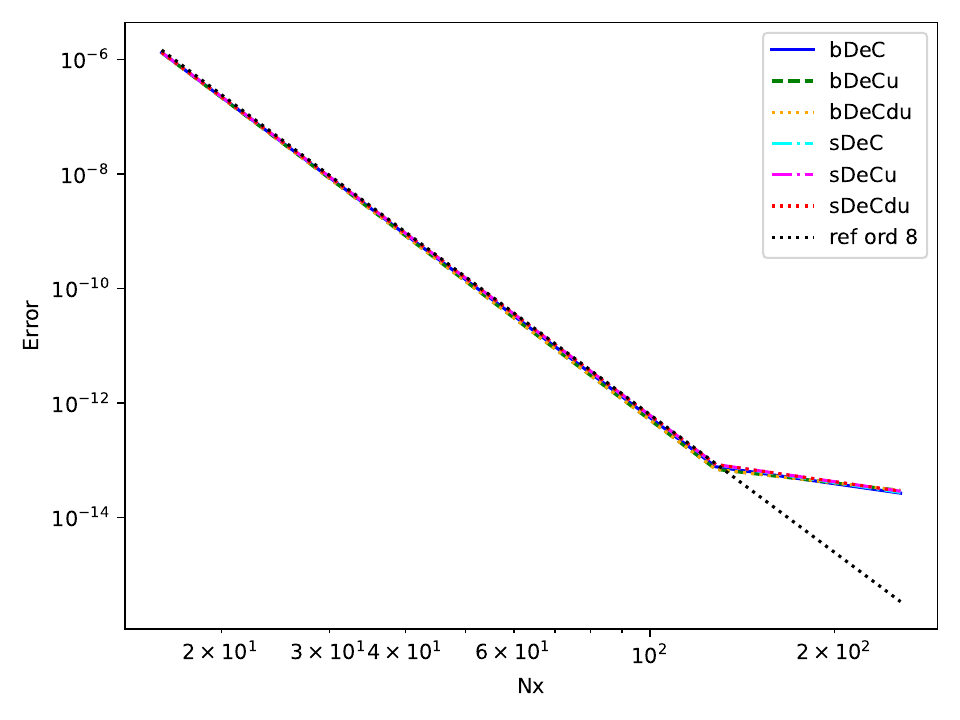}
	\includegraphics[width=0.49\textwidth]{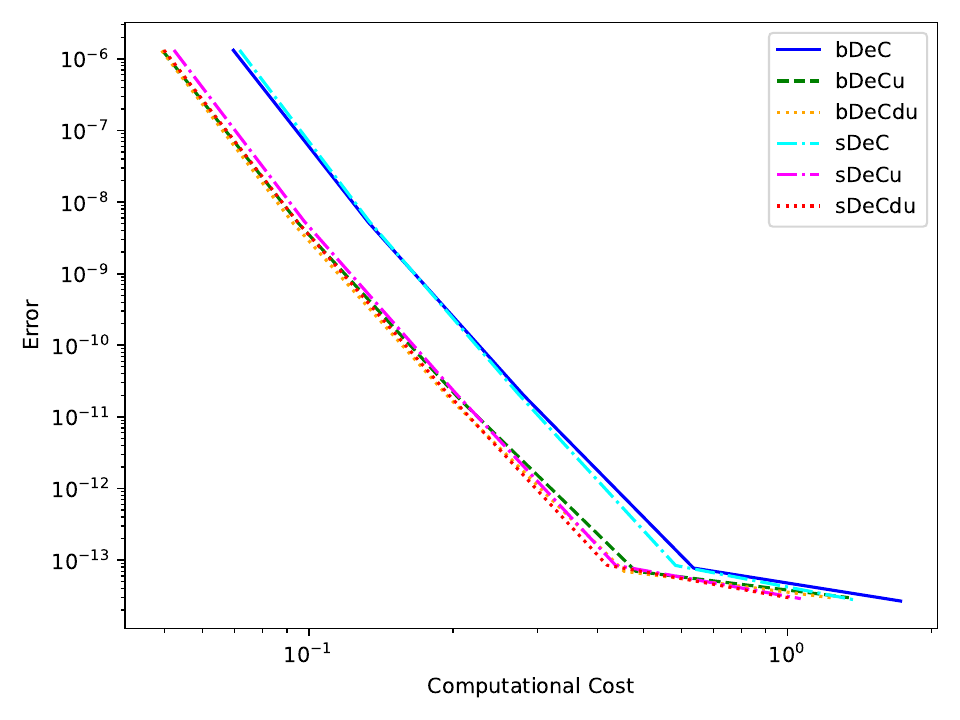}
	\caption{Advection--diffusion equation: Convergence and efficiency analysis of the IMEX DeC schemes. Left: Convergence; Right: efficiency. Top: order 5; Bottom: order 8.}\label{fig:adv_diff}
\end{figure}

We test the convergence of the methods for the advection-diffusion problem with $a=0.2$, $\mu=0.01$, and periodic boundary conditions on the domain $[0,1]$. The initial condition is set to $u(x,0)=\sin(2\pi x)$, with corresponding exact solution $u(x,t)=\exp(-4\pi^2\mu t)\sin\bigl(2\pi(x-at)\bigr).$
We run the simulations until the final time $T_f=1$, and we use a hyperbolic CFL condition to determine the time-step size, with $\Delta t = \text{CFL} \cdot \Delta x / \vert a \vert $, where $\text{CFL}=1$.

In Figure~\ref{fig:adv_diff}, we report the errors of the methods against time-step size and computational time for orders 5 and 8. The results show that the expected order of accuracy is achieved for all methods, with the efficient modifications being faster than the original methods, in particular for high order. For a fixed mesh discretization, the bDeCu and bDeCdu methods are again more accurate than the bDeC method, while the errors of the sDeCu and sDeCdu methods are very close to the sDeC ones.
In all cases, however, we have a clear computational advantage as can be inferred from the efficiency plots.

\subsection{Allen-Cahn}\label{sec:allen_cahn}
In this section, we consider the two--dimensional Allen-Cahn equation, introduced in~\cite{allen1979microscopic}, reading
\begin{equation}
	u_t+\uvec{a}\cdot \nabla_{\uvec{x}} u=\mu \Delta_{\uvec{x}}u+r(u-u^3), \quad (\uvec{x},t)\in \Omega\times[0,T_f],
	\label{eq:AC}
\end{equation}
where $\mu,r\geq 0$ are the diffusion and reaction coefficients, respectively, while $\uvec{a}\in \mathbb{R}^2$ represents a constant background advection field.
In this case, indeed, $\uvec{x}\in\Omega\subset \mathbb{R}^2$. 
This model is mainly used in the context of phase separation and interface motion in materials.
More specifically, $u$ is an order parameter taking values between -1 and +1, representing two different phases of a material. 
The reaction term, $r(u-u^3)$, pushes the solution towards those two stable states; while the diffusion term, $\mu \Delta_{\uvec{x}}u$, smoothens the transition layer between them.

Herein, we set $\mu:=0.01$, $r:=10$, $\uvec{a}:=(0.2,0.1)^\top$.
For the test, we consider the computational space domain $[0,2\pi]\times [0,\pi]$ with periodic boundary conditions, we prescribe a random initial condition uniformly distributed on $[-0.05,0.05]$, and run our simulations until the final time $T_f:=1$.
The spatial discretization is performed through a Cartesian extension, dimension by dimension, of the structures defined for the previous test.
In particular, the univariate spatial derivative operators described above are also applied dimension by dimension, notice that both components of $\uvec{a}$ are positive. 
Just like before, the advection part is treated explicitly. The whole right-hand side of~\eqref{eq:AC}, consisting of the diffusion and of the nonlinear reaction term, is treated implicitly. The time-step is chosen as $\Delta t=0.022$ tuned so that the nonlinear solver converges in a reasonable number of iterations.
Of course, it also satisfies the hyperbolic CFL condition: 
\begin{equation}
	\Delta t
	\leq 
	\min\left(
	\frac{\Delta x}{|a_x|},
	\frac{\Delta y}{|a_y|}
	\right).
\end{equation}

\begin{figure}
	\begin{minipage}{0.49\textwidth}
		\centering
		sDeC order 2
	\end{minipage}
	\begin{minipage}{0.49\textwidth}
		\centering
		sDeCu order 3
	\end{minipage}\\
	\includegraphics[width=0.49\textwidth]{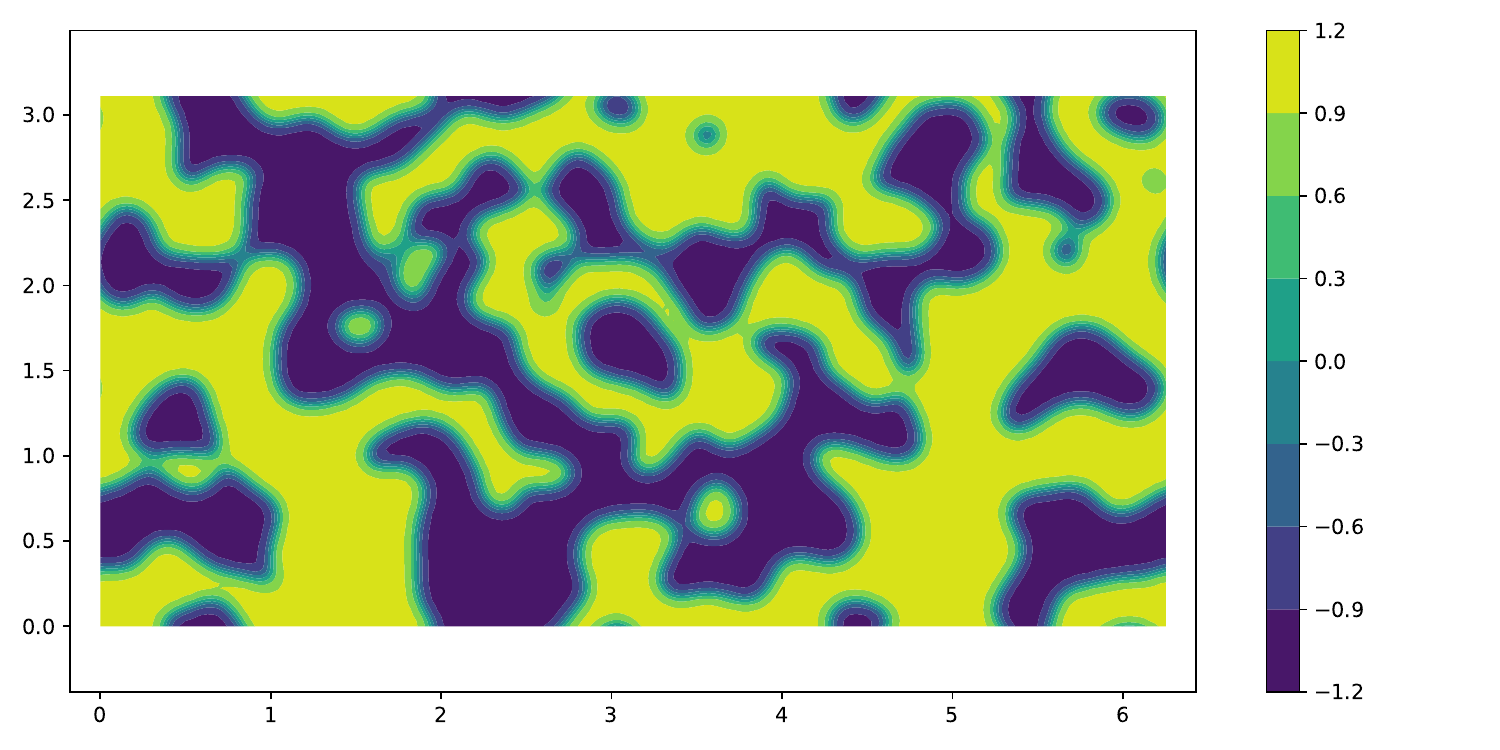}
	\includegraphics[width=0.49\textwidth]{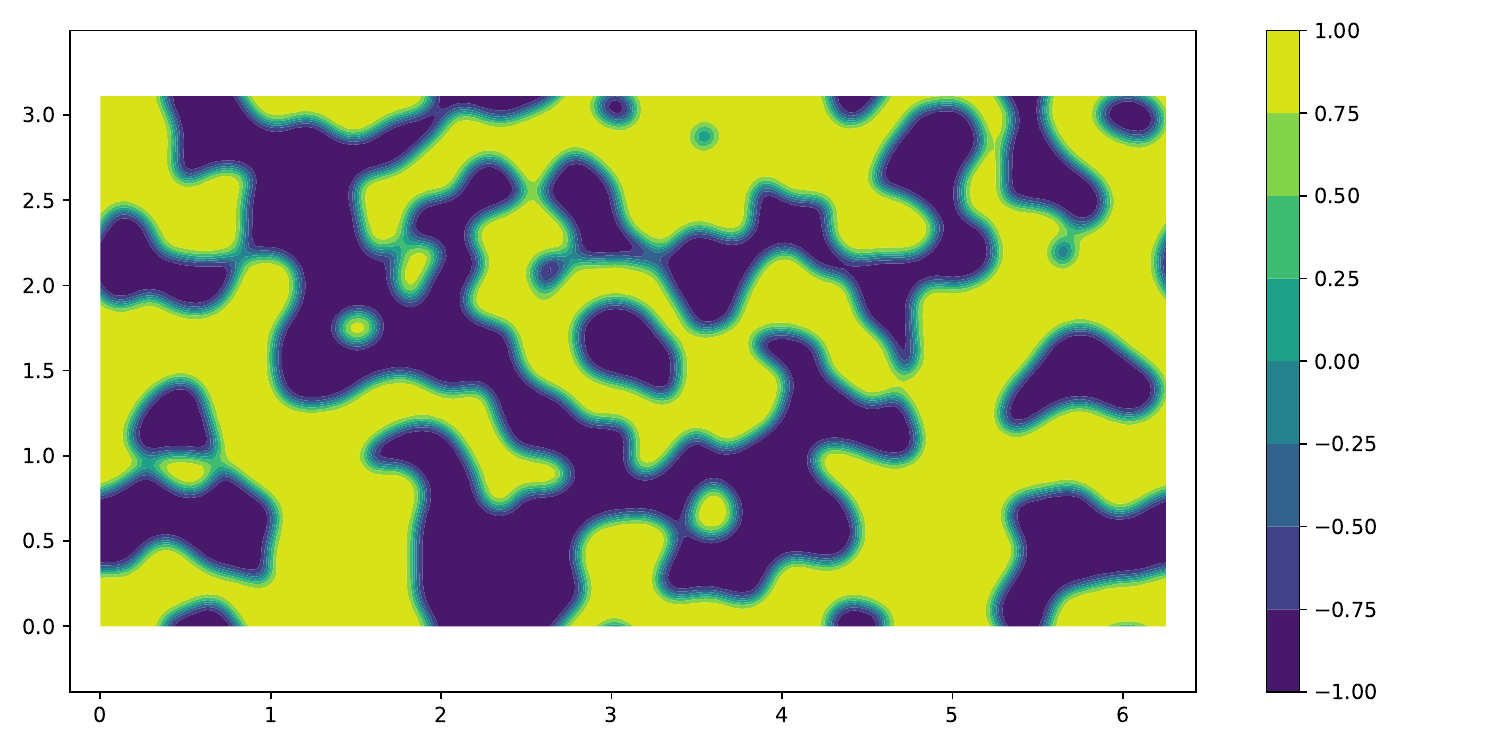}\\
		\begin{minipage}{0.49\textwidth}
		\centering
		sDeCu order 5
	\end{minipage}
	\begin{minipage}{0.49\textwidth}
		\centering
		Adaptive sDeCu order 5
	\end{minipage}\\
	\includegraphics[width=0.49\textwidth]{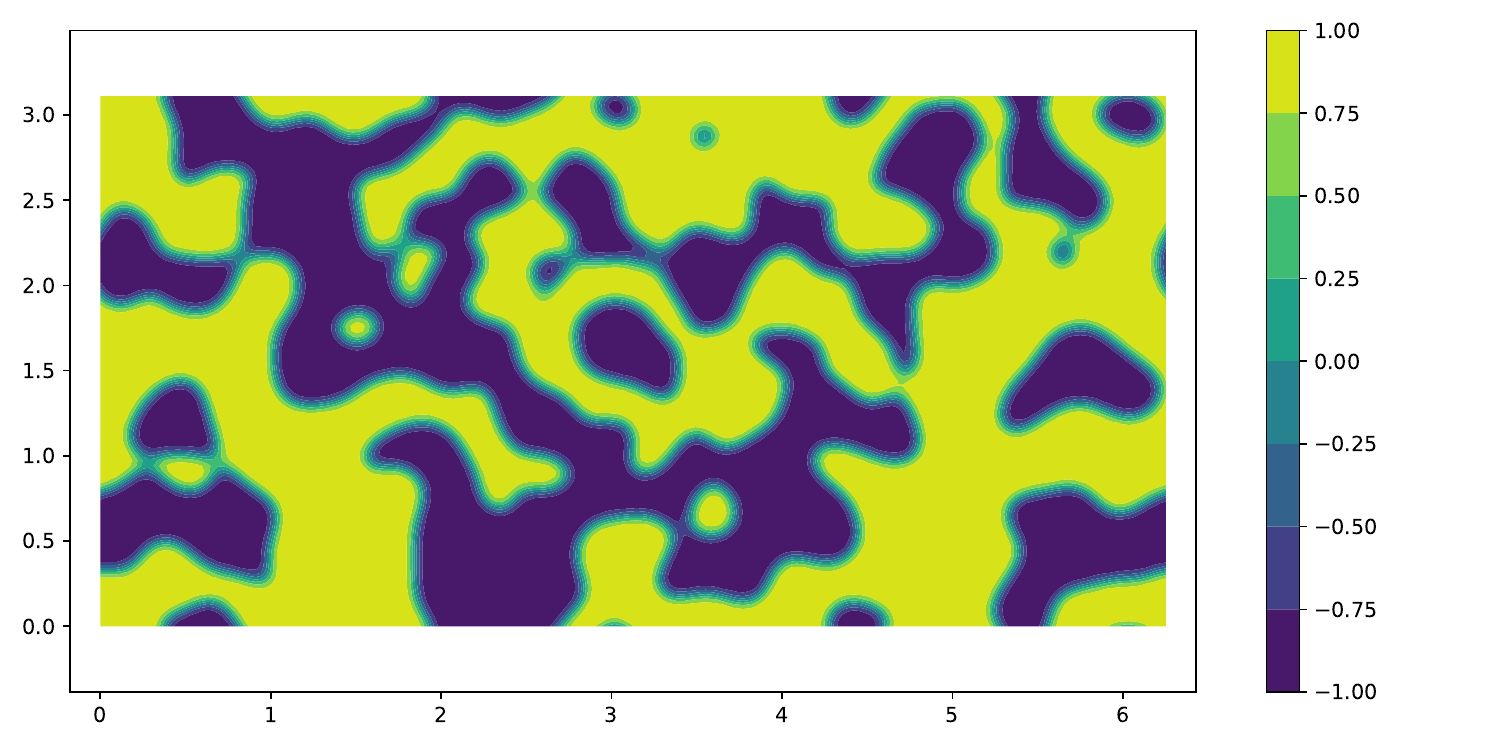}
	\includegraphics[width=0.49\textwidth]{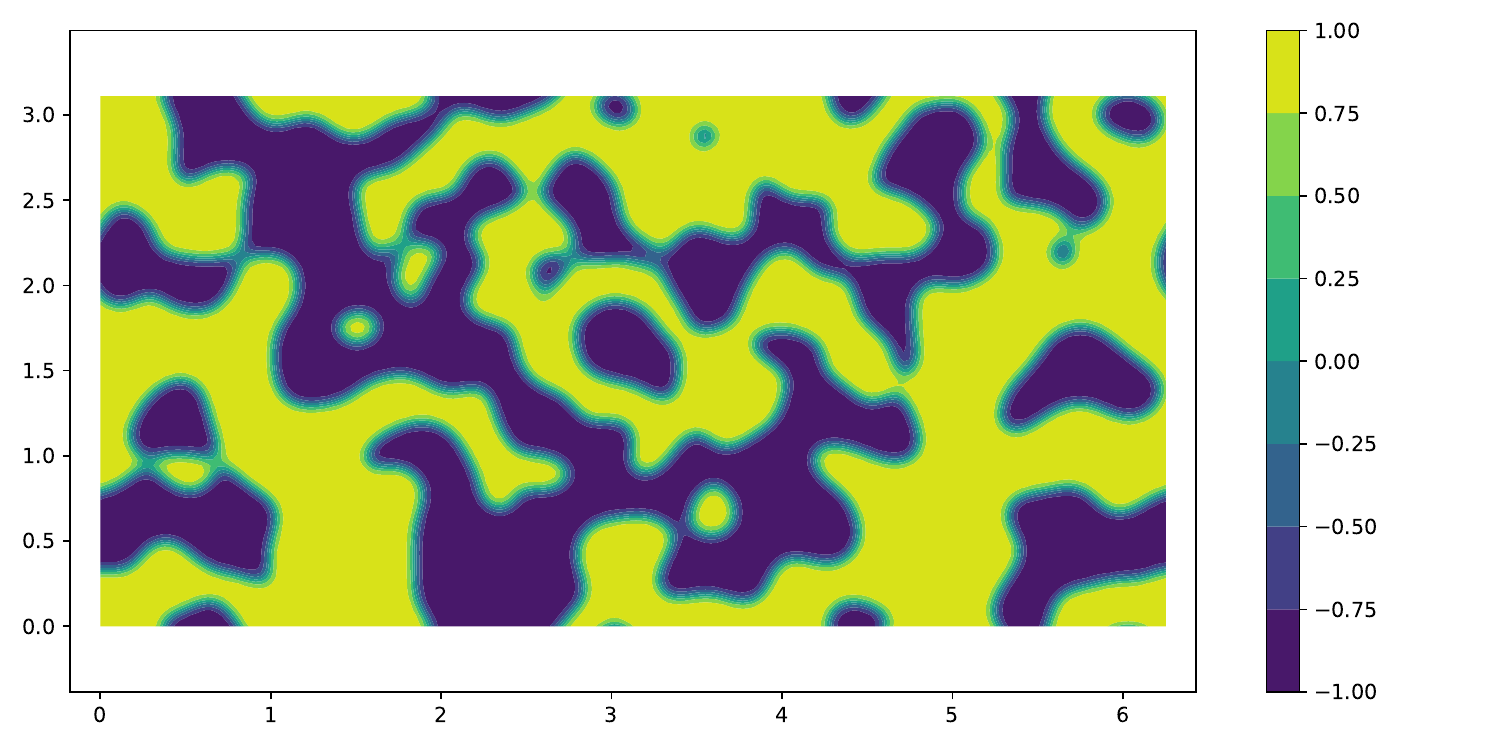}
	\caption{Allen--Cahn: Results obtained on a grid of $200\times 100$ points. The adaptive version makes use of fixed fifth-order of accuracy in space and $\varepsilon:=10^{-3}$.} \label{fig:AC}
\end{figure}

In Figure~\ref{fig:AC}, we report the solution obtained with the second-order sDeC scheme, with the third- and fifth-order sDeCu schemes and with the adaptive sDeCu scheme of (maximum) order 5, where the adaptation is only performed in time, according to Equation~\eqref{eq:tol} with $\varepsilon:=10^{-3}$, using a fixed spatial order of 5. Indeed, the final number of subtimenodes corresponds to order 5, but fewer iterations are performed if the convergence tolerance is matched.

\begin{figure}
	\centering
	\begin{tikzpicture}
		
		\begin{axis}[
			width=0.8\textwidth,
			height=0.45\textwidth,
			xlabel={Time},
			ylabel={Number of iterations},
			xmin=0,
			grid=major,
			legend style={
				at={(0.98,0.98)},
				anchor=north east,
				draw=none,
				fill=none
			},
			]
			
			\addplot[
			thick,
			blue,
			mark=x,
			mark size=1.2pt,
			]
			table[
			x index=0,
			y index=1,
			col sep=space,
			header=true,
			]{./figures/AllenCahn2D_DeC_small_sub_staggered_adaptive_Nx_200_p_times.txt};
			
			
		\end{axis}
	\end{tikzpicture}
	
	\caption{Allen--Cahn: Number of DeC iterations performed by the adaptive sDeCu method as a function of time with error tolerance $\varepsilon=10^{-3}$.}
	\label{fig:adaptive_iterations_allen_cahn}
\end{figure}
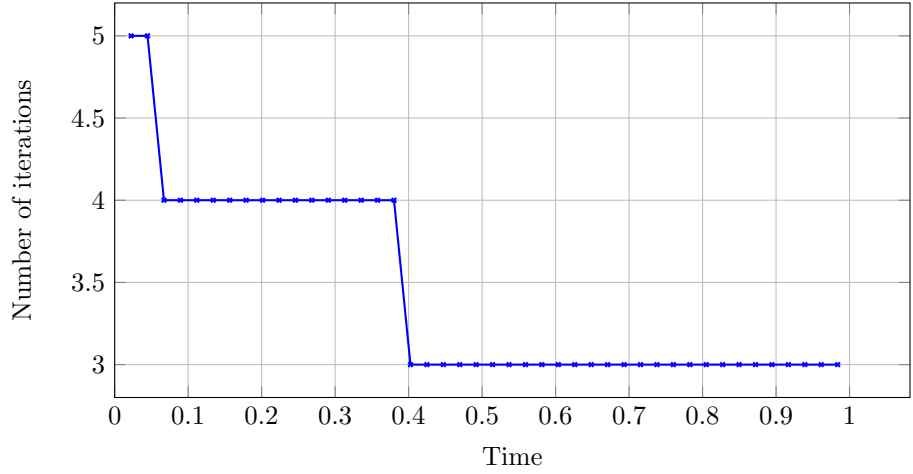

It is interesting to see that all methods converge to the same solution even if starting from the same random initial condition. 
The second-order method is clearly less accurate than the higher-order ones, as we can notice from the less sharp interfaces between the two phases. 
The difference between the third- and the fifth-order methods is qualitatively negligible. 
The adaptive method is able to reach the same qualitative accuracy as the fifth-order method, while using less iterations in time, hence being more efficient.
In Figure~\ref{fig:adaptive_iterations_allen_cahn}, we report the number of iterations used in each time-step for the adaptive method. The number of iterations is compliant with the dynamics of the test: it is higher at the beginning of the simulation, when the solution is still far from the two stable states, and it quickly decreases as the solution approaches an equilibrium between the phases, while the advection operator always moves the solution.

The computational cost for the fifth-order sDeC method was 253 seconds, for sDeCdu was 210 seconds, while for the adaptive method it was 143 seconds, saving around 77\% of the computational time with respect to the original version and 30\% of the computational time with respect to sDeC.

\subsection{Cahn-Hilliard}\label{sec:cahn_hilliard}
In this section, we consider the two--dimensional Cahn--Hilliard equation, introduced in~\cite{cahn1958free}, reading
\begin{equation}
	u_t+\uvec{a}\cdot \nabla_{\uvec{x}} u=\mu \Delta_{\uvec{x}}(u^3-u-\gamma \Delta_{\uvec{x}} u), \quad (\uvec{x},t)\in \Omega\times[0,T_f].
\end{equation}
In this equation, $\mu,\gamma\geq 0$ represent the mobility and interface-energy coefficients, and $\uvec{a}\in\mathbb{R}^2$ is again a constant background advection field.
This model, pre-existing with respect to Allen-Cahn, is used to describe phase separation in binary mixtures through a diffuse-interface formulation. In contrast with the previous model, the Cahn-Hilliard equation is mass-conservative.
Again, $u$ is an order parameter whose values lie between $-1$ and $+1$, representing the two different phases of the material.
The quantity $u^3-u-\gamma \Delta_{\uvec{x}}u$ represents the chemical potential.
The nonlinear term, $u^3-u$, energetically favors the two phases, while the term $-\gamma \Delta_{\uvec{x}}u$ smoothens the transition layer between them.
The outer Laplacian drives the redistribution of the order parameter according to spatial variations of the chemical potential, thereby preserving its total mass.

In this case, we set $\mu:=1$, $\gamma:=0.001$, $\uvec{a}:=(0.2,0.1)^\top$, and we consider the same domain, initial and boundary conditions as for the previous test, with final time $T_f:=2$.
The spatial discretization is, like in the previous test, obtained by applying the one-dimensional derivative operators dimension by dimension. In particular, the biharmonic operator is obtained by applying the discrete Laplacian twice.
%
As for the Allen--Cahn problem, the advection is treated explicitly, while, the other spatial terms are handled implicitly. Furthermore, we use the same fixed time-step $\Delta t = 0.022$  for the computations.

\begin{figure}
	\centering
	\includegraphics[width=0.49\textwidth]{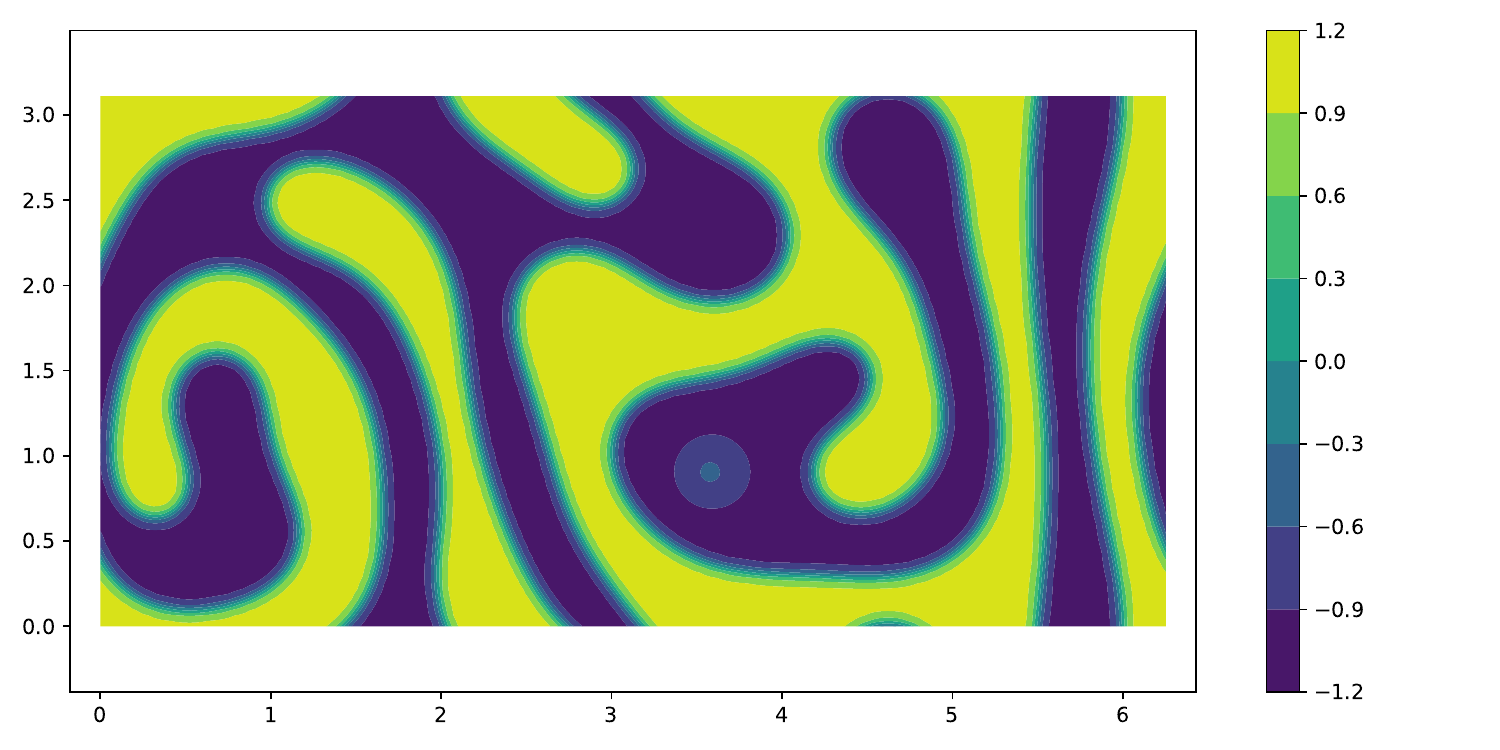}
	\begin{tikzpicture}
		
		\begin{axis}[
			width=0.44\textwidth,
			height=0.25\textwidth,
			xlabel={Time},
			ylabel={Number of iterations},
			xmin=0,
			grid=major,
			]
			\addplot[
			thick,
			red,
			dashed,
			mark=o,
			mark size=1.2pt,
			]
			table[
			x index=0,
			y index=1,
			col sep=space,
			header=true,
			]{./figures/CahnHilliard2D_DeC_small_sub_staggered_adaptive_Nx_200_p_times.txt};
			
		\end{axis}
	\end{tikzpicture}
	\caption{Cahn-Hilliard problem. Results at final time $T_f=2$ obtained on a grid of $200\times 100$ points with adaptive sDeCu of order 5 for $\varepsilon:=10^{-3}$, employing fixed fifth-order of accuracy in space (left). Number of DeC iterations performed by the adaptive sDeCu method as a function of time with error tolerance $\varepsilon=10^{-3}$ (right).} \label{fig:CH}
\end{figure}
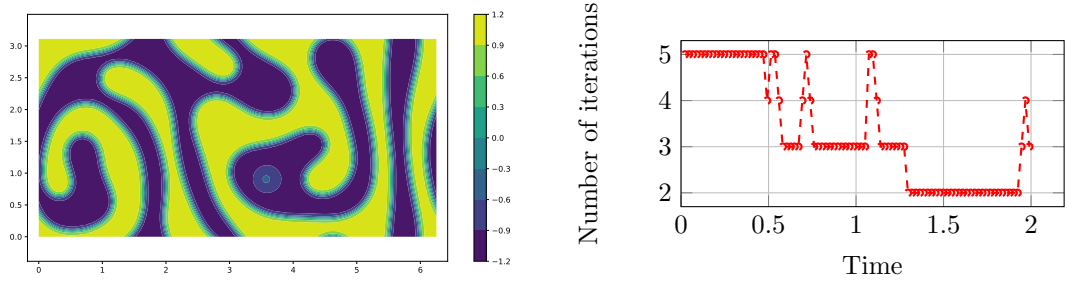

In this case, we focus on the adaptive version of the algorithm. More specifically, we consider the adaptive sDeCu of the previous test with (maximum) time order 5, with $\varepsilon:=10^{-3}$, and constant space order 5.
The obtained final solution, for a mesh with 200$\times$100 mesh nodes, is displayed in Figure~\ref{fig:CH} (left), featuring patterns consistent with the ones reported in existing literature.

\begin{figure}
	\begin{minipage}{0.32\textwidth}
		\centering
		Time $=$ 0.693\\
		\includegraphics[width=\textwidth, trim={70 20 80 45},clip]{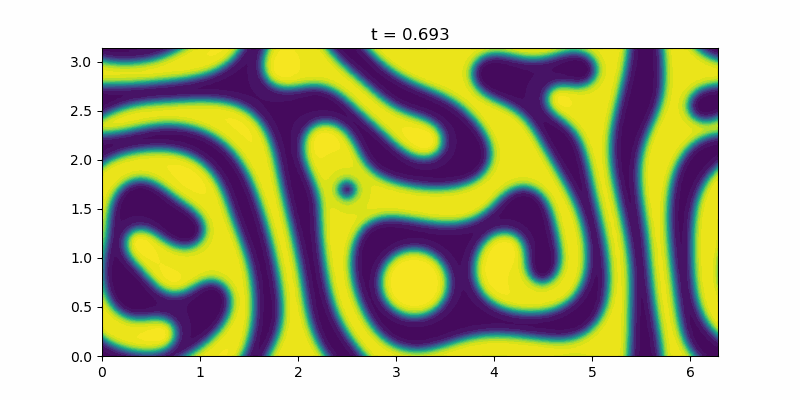}\\
		Time $=$ 0.716\\
		\includegraphics[width=\textwidth, trim={70 20 80 45},clip]{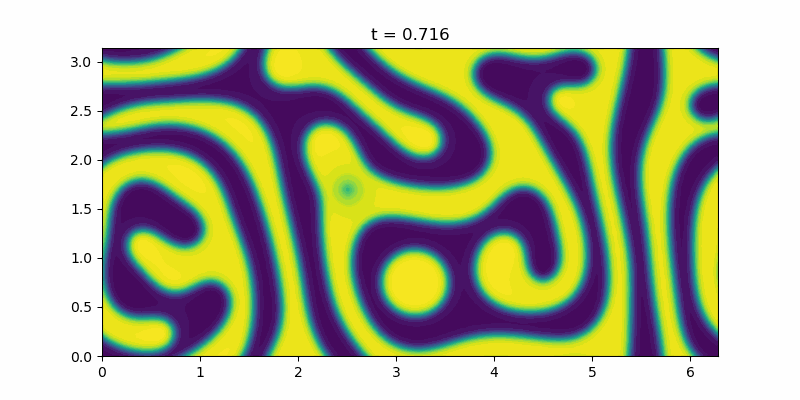}\\
		Time $=$ 0.738\\
		\includegraphics[width=\textwidth, trim={70 20 80 45},clip]{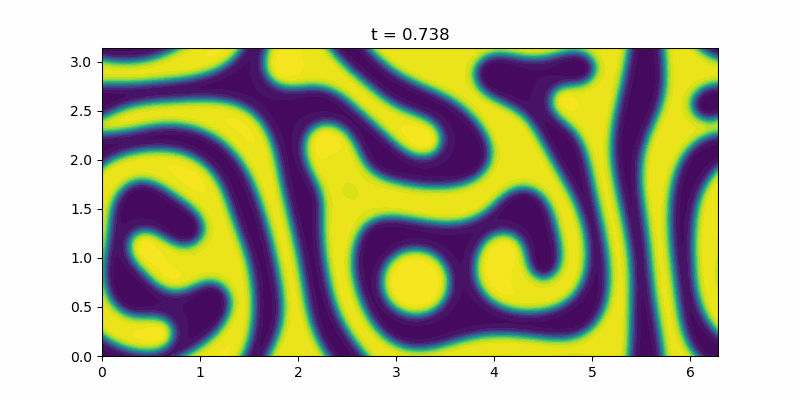}
	\end{minipage}
	\begin{minipage}{0.32\textwidth}
		\centering
		Time $=$ 1.073\\
		\includegraphics[width=\textwidth, trim={70 20 80 45},clip]{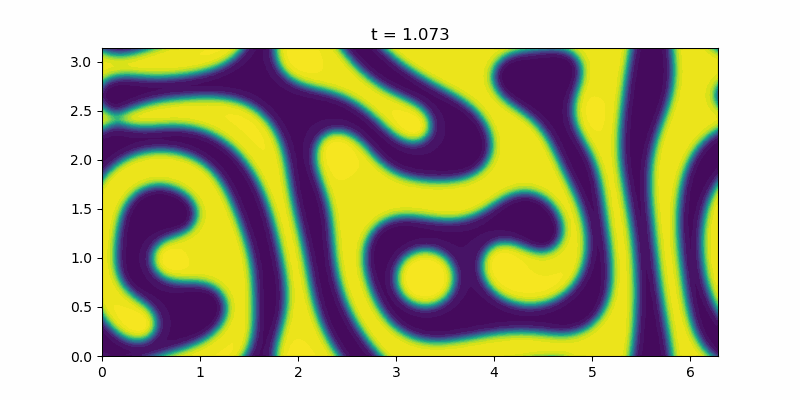}
		Time $=$ 1.096\\
		\includegraphics[width=\textwidth, trim={70 20 80 45},clip]{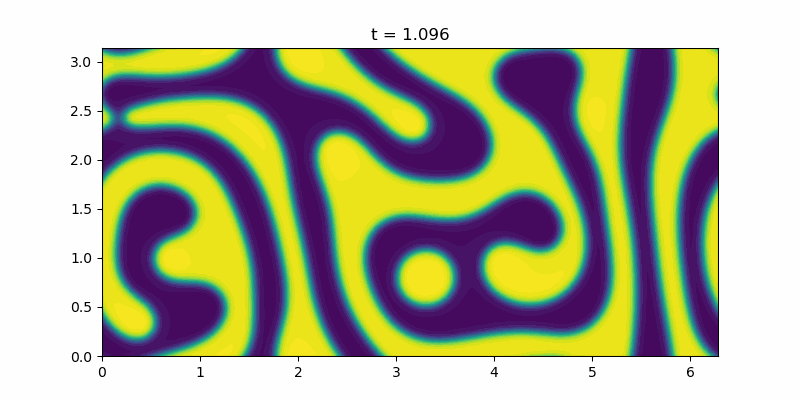}
		Time $=$ 1.118\\
		\includegraphics[width=\textwidth, trim={70 20 80 45},clip]{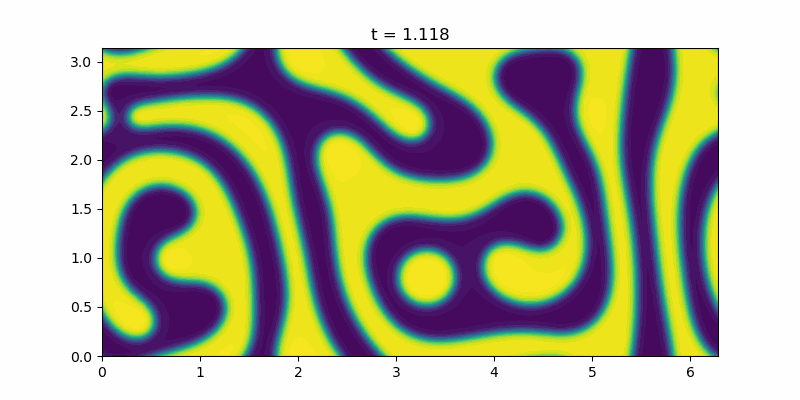}		
	\end{minipage}
	\begin{minipage}{0.32\textwidth}
		\centering
		Time $=$ 1.945\\
		\includegraphics[width=\textwidth, trim={70 20 80 45},clip]{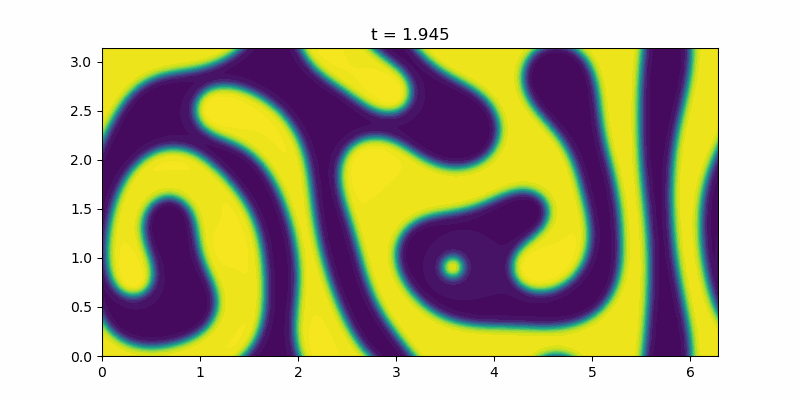}
		Time $=$ 1.969\\
		\includegraphics[width=\textwidth, trim={70 20 80 45},clip]{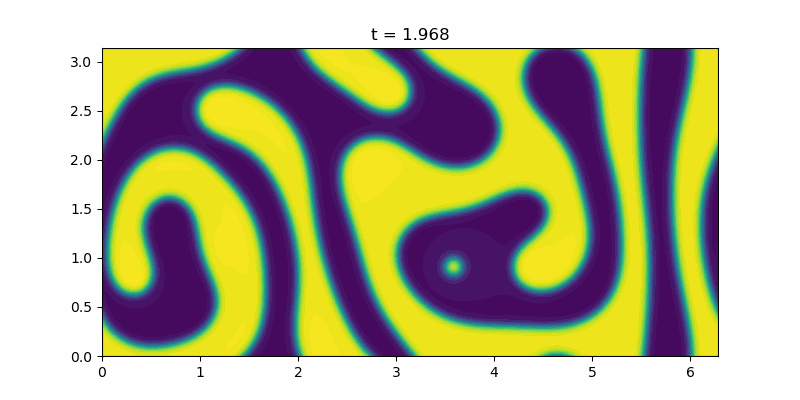}
		Time $=$ 1.990\\
		\includegraphics[width=\textwidth, trim={70 20 80 45},clip]{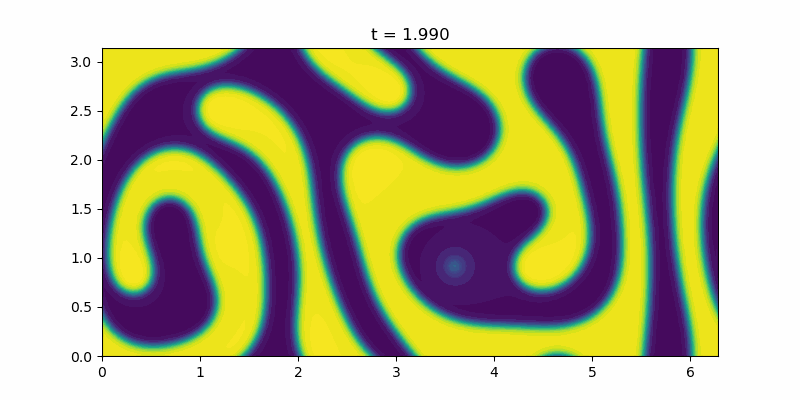}				
	\end{minipage}
	\caption{Cahn--Hilliard: Different snapshots of the adaptive sDeCu method corresponding to local increases in time of the number of iterations, i.e., of the order of accuracy. These correspond to the disappearance of a blue point close to the coordinates (2.5,1.7) at $t\sim 0.7$ (left column), the detachment of two yellow branches close to the coordinates (0.2,2.5) at $t\sim 1.1$ (central column) and the disappearance of a yellow area close to the coordinates (3.5,0.9) at $t\sim 1.95$ (right column).}
	\label{fig:snapshots}
\end{figure}

In Figure~\ref{fig:CH} (right), we report the number of iterations used in each time-step for the adaptive strategy. The number of iterations varies along the simulation according to the dynamics of the solution. In particular, the approach is able to detect when to use more or fewer iterations, hence being more efficient.
Some snapshots of the solution at three different times, corresponding to an increase of the number of DeC iterations, are reported in Figure~\ref{fig:snapshots}, where we can appreciate how they all correspond to key pattern changes.
Finally, we remark that while the original sDeC of order 5 took 3734 seconds to run the simulation, the sDeCu method required 2518 seconds and the adaptive sDeCu method of order 5 required 2006 seconds. Hence, the adaptive strategy brought a saving of 20\% of the computational cost with respect to the sDeCu and of 86\% with respect to the original sDeC method.
This test shows the potential of the proposed adaptive strategy. Further applications, where adaptivity is performed in space and in time, are left for future works.


\section{Conclusions}\label{sec:conclusions}

In this work, we have investigated efficient modifications of two families of implicit--explicit (IMEX) Deferred Correction (DeC) methods of arbitrary high order.
The proposed modifications are based on introducing interpolation processes between consecutive DeC iterations, so that the order of the discretization structures employed at each iteration matches the accuracy attained at that stage.
This allows the first low-order iterations to be performed using cheaper structures while progressively increasing the approximation order throughout the iterative process, leading to more efficient schemes.
%
%
We have considered both solution-based and right-hand-side-based interpolation strategies and assessed their accuracy, stability and computational efficiency.

The stability analysis shows that the proposed modifications preserve, to a large extent, the stability properties of the corresponding original IMEX DeC schemes.
The numerical experiments confirm that the expected orders of accuracy are achieved and that the modified methods provide a clear computational advantage over the original formulations, especially at high orders.
This behavior has been observed for problems of different nature, including moderately stiff and highly stiff vibrating systems, the rescaled Van der Pol oscillator, and semidiscretizations of advection--diffusion, Allen--Cahn and Cahn--Hilliard partial differential equations (PDEs).

The iterative structure of the proposed schemes has also been exploited to construct adaptive methods in which the number of iterations, and therefore the temporal order of accuracy, is selected according to a prescribed tolerance.
The numerical results show that these adaptive formulations are able to adjust the number of iterations to the dynamics of the solution and to achieve accuracy comparable to that of fixed high-order methods at a reduced computational cost.
This behavior is particularly evident in the Allen--Cahn and Cahn--Hilliard tests, where the adaptive strategy performs more iterations during the most active phases of the evolution and fewer iterations when the dynamics become smoother.

Further developments will concern more general adaptive criteria also involving spatial adaptivity, see~\cite{micalizzi2023efficient,micalizzi2025smoothness,micalizzi2026adaptive} and references therein, and applications to more involved multiscale PDEs in the context of asymptotic-preserving schemes, see for example~\cite{chertock2026new,chertock2026asymptotic}.

\section*{Declarations}

\subsection*{Funding}
L.M. is funded by the LeRoy B. Martin, Jr. Distinguished Professorship Foundation.
D.T. is a member of the INdAM GNCS group in Italy and was supported by the Ateneo Sapienza project 2024 ``Advanced Computational Methods for Real-World Applications: Data-Driven Models, Hyperbolic Equations and Optimal
Control''.
The authors gratefully acknowledge Sapienza University of Rome for financial support through the Visiting Professor Programme of the Department of Mathematics ``Guido Castelnuovo'', funded within the Progetto di Eccellenza 2023--2027 (project code 282933\_DIP\_ECC\_2023\_2027\_029, CUP B83C23001390001), and through the Sapienza Visiting Professors Programme 2025 for joint research activities.
A substantial part of the work was developed during the visit of L.M. to the Department of Mathematics.

\subsection*{Competing interests}
The authors have no relevant financial or non-financial interests to disclose.

\subsection*{Data availability}
The data generated during the current study are available from the corresponding author upon reasonable request.

\bibliography{sn-bibliography}
\bibliographystyle{plain}

\end{document}